\documentclass[11pt,twoside]{article}
\usepackage{mathrsfs}
\usepackage{bm}
\usepackage{stmaryrd}
\usepackage{amsfonts,amsmath}
\begin{document}

\title{\bf Rayleigh-Taylor instability for compressible\\
 rotating flows }

\author{{ Ran Duan$^a$, \ Fei Jiang$^b$, \ Song Jiang$^b$ }\\[2mm]
\small$^a$ School of Mathematics and Statistics,\\
\small Central China Normal University, Wuhan 430079, China.\\[2mm]
\small $^b$ Institute of Applied Physics and Computational
\small Mathematics,\\ \small Beijing 100088, China.}

\date{}
\vskip 0.2cm

\maketitle
\renewcommand{\thefootnote}{\fnsymbol{footnote}}

\footnotetext[2]{The research of Ran Duan was supported by three
grants from the National Natural Science Foundation of China (Grant
No. 11001096) , the Special Fund for Basic Scientific Research  of
Central Colleges (Grant No. CCNU10A01021), and  China Postdoctoral
Science Foundation (Grant No.20110490327); the research of Fei Jiang
was supported by the Fujian Provincial Department of Science and
Technology (Grant No. JK2009045) and NSFC (Grant No. 11101044); the
research of Song Jiang was supported by the National Basic Research
Program under the Grant 2011CB309705 and NSFC (Grant No.40890154)}
\footnotetext[3]{Email addresses: duanran@mail.ccnu.edu.cn (Ran Duan),
jiangfei0591@163.com (Fei Jiang), jiang@iapcm.ac.cn (Song Jiang)}

\vskip 0.2cm \arraycolsep1.5pt
\newtheorem{Lemma}{Lemma}[section]
\newtheorem{Theorem}{Theorem}[section]
\newtheorem{Definition}{Definition}[section]
\newtheorem{Proposition}{Proposition}[section]
\newtheorem{Remark}{Remark}[section]
\newtheorem{Corollary}{Corollary}[section]

\begin{abstract}

In this paper, we investigate the Rayleigh-Taylor instability problem for two
compressible, immiscible, inviscid flows rotating with an constant
angular velocity, and evolving with a free interface in the presence
of a uniform gravitational field. First we construct the
Rayleigh-Taylor steady-state solutions with a denser fluid lying
above the free interface with the second fluid, then we turn to an
analysis of the equations obtained from linearization around such a
steady state. In the presence of uniform rotation, there is no
natural variational framework for constructing growing mode
solutions to the linearized problem. Using the general method of
studying a family of modified variational problems introduced in
\cite{Y-I2}, we construct normal mode solutions that grow
exponentially in time with rate like $e^{t\sqrt{c|\xi|-1}}$, where
$\xi$ is the spatial frequency of the normal mode and the constant
$c$ depends on some physical parameters of the two layer fluids. A Fourier
synthesis of these normal mode solutions allows us to construct
solutions that grow arbitrarily quickly in the Sobolev space $H^k$, and
lead to an ill-posedness result for the linearized problem. Moreover, from the analysis
we see that rotation diminishes the growth of instability.
Using the pathological solutions, we then demonstrate the
ill-posedness for the original non-linear problem in some sense.

\bigbreak


\bigbreak

\end{abstract}

\section{Introduction}
\setcounter{equation}{0}

 Hydrodynamic instabilities at the interface of two materials
of different densities are a critical issue in high energy density
physics (HEDP). The Rayleigh-Taylor instability (RTI) occurs when a
fluid accelerates another fluid of high density \cite{L.R1,L.R,G.I}.
The RTI is ubiquitous in HEDP, such as high Mach number shocks and
jets, radiative blast waves and radioactively driven molecular
clouds, gamma-ray bursts and accreting black holes, etc. Due to the
importance in physics mentioned above, there have been many studies
related to RTI from both physical and numerical simulation points of
view in the literature, but only few analytical results.
One classical case of RTI is to consider two
completely plane-parallel layers of immiscible inviscid fluids, the heavier
on top of the light one and both subject to the earth's gravity. The
equilibrium here is unstable to certain perturbations or
disturbances. For this model, Y. Guo and I. Tice \cite{Y-I}
established a variational framework for nonlinear instabilities.
They constructed solutions that grow arbitrarily quickly in the
Sobolev space by the method of Fourier synthesis, which leads to an
ill-posedness to the perturbed problem. Ever since then, several
works studied the influences of different physical quantities to the
linear RTI, such as the stabilized effect of viscosity and surface
tension  (see \cite{Y-I2,FJSJWW}), and the effect of magnetic field
(see \cite {D-J-J,FJSJYW}), etc.

 While, if the two fluids are all subject to a uniform rotation with an constant
angular velocity, how does this rotation influent the RTI?
Stabilize it or make it more unstable? To the best of our
knowledge, physicists \cite{S. Chandrasekhar} showed that rotation
can only slow down the growth rate of the disturbance a little bit, but not
prevent the incompressible fluids from becoming unstable. Now, an interesting question is
whether it is still true for compressible flows?
In this paper, for inviscid compressible flows without the centrifugal force, we first
prove that the linearized system is unstable in the Hadamard sense. Then,
we establish the ill-posedness for the original non-linear
problem in some sense.

Next, we formulate the problem in details for further discussion.
%%%%%%%%%%%%%%%%%%%%%%%%%%%%%%%%%
\subsection{Formulation in Eulerian coordinates}

We suppose that the fluids are confined between two rigid planes.
As in \cite{Y-I}, we denote this
infinite slab by ${\Omega}={\mathbb{R}}^{2}\times (-m,l)\subset
{\mathbb{R}}^3$. The fluids are
separated by a moving free interface $\Sigma (t)$ ($t\geq 0$) that extends to
infinity in every horizontal direction. The interface divides
$\Omega$ into two time-dependent, disjoint, open subsets
$\Omega_{\pm}(t)$, so that $\Omega=\Omega_+(t)\cup \Omega_-(t)\cup
\Sigma (t)$ and $\Sigma (t)=\bar{\Omega}_+(t)\cap
\bar{\Omega}_-(t)$. The motion of the fluids is driven by the
constant gravitational field along $e_3$-the $x_3$ direction,
$\mathbf{g}=(0,0,-g)$ with $g>0$ and the rotation with an angular velocity
$\boldsymbol{\omega}=(0,0,\omega)$ about the vertical direction. Quantities
describing fluids are their density and velocity, which are given for each
$t\geq 0$ by, respectively,
\begin{eqnarray*}  % \label{0101}
\rho_\pm(t,\cdot):\Omega_\pm(t)\rightarrow {\mathbb{R}^+}, \quad
u_\pm(t,\cdot):\Omega_\pm(t)\rightarrow {\mathbb{R}^3}.
\end{eqnarray*}
We shall assume that at a given time $t\geq 0$ these functions have
well-defined traces onto $\Sigma(t)$.

 The fluids are governed by the following equations:
\begin{equation}\label{0102}
\left\{
\begin{array}{ll}\partial_t\rho_{\pm}+ \mathrm{div}(\rho_\pm u_\pm)=0, \\[1mm]
\rho_\pm(\partial_t u_{\pm}+u_{\pm}\cdot \nabla
u_{\pm})+2\rho_{\pm}(\boldsymbol{\omega}\times u_{\pm})+\nabla
p(\rho_{\pm})=-g\rho_{\pm}e_3,
                  \end{array}
                \right.
\end{equation} for $t>0$ and $x\in \Omega_\pm(t)$.
Here we have written $g>0$ for the gravitational constant,
$e_3=(0,0,1)$ for the vertical unit vector, and $-g e_3$ for the
acceleration due to gravity, $2\rho_\pm(\boldsymbol{\omega}\times u_{\pm})$
represents the Coriolis force, while the centrifugal force
$\rho_\pm\nabla|\boldsymbol{\omega}\times x|^2/2$, like in
\cite{AKSSTNS,Feireisl}, is neglected on the basis of a standard
argument which indicates that the motion is dominated by gravitation and
the centrifugal force is small.

We assume a general barotropic pressure law of the form
$p_\pm=p_\pm(\rho)\geq 0$ with $p_\pm\in C^\infty(0,\infty)$ and
strictly increasing. We will also assume that $1/p'_\pm\in
L^\infty_{{\rm{loc}}}(0,\infty)$. Finally, in order to construct a
steady state solution with an upper fluid of greater density at
$\Sigma(t)$, we will assume that
\begin{eqnarray*}  %%\label{001}
Z:=\left\{z\in(0,\infty)~|~ p_-(z)>p_+(z) \ {\rm{and}}\  p_-(z)\in p_+(0,\infty)\right\}
\neq \emptyset.   \end{eqnarray*}
In particular, this requires that the pressure laws be distinct, i.e. $p_-\neq p_+$.

Now, we prescribe the jump conditions that, from the physical point of view, both normal component of
the velocity and the pressure are continuous across the free interface (since no surface tension
is taken into account), see \cite{S. Chandrasekhar, J.Wehausen-1960}.
Therefore, the jump conditions at the free interface read as
$$ [u\cdot \nu]|_{\Sigma (t)}=0, \qquad [p]|_{\Sigma (t)}=0  $$
for each $t>0$, where we have denoted the normal vector to $\Sigma(t)$ by $\nu$, the
trace of a quantity $f$ on $\Sigma(t)$ by $f|_{\Sigma(t)}$ and the interfacial jump by
\begin{eqnarray*}   %%%\label{0106}
[f]|_{\Sigma(t)}:=f_+|_{\Sigma(t)}-f_-|_{\Sigma(t)}.
\end{eqnarray*}

On the fixed boundaries, we consider that the normal component of the
fluid velocity vanishes, that is,
\begin{eqnarray*} % \label{0107}
u_+(t,x',-m)\cdot e_3=u_-(t,x',l)\cdot e_3=0,\mbox{ for all }
t\geq 0,\ x':=(x_1,x_2)\in \mathbb{R}^{2}.
\end{eqnarray*}

The motion of the free interface is coupled to the evolution
equations for the fluids (\ref{0102}) by requiring that the interface
be advected with the fluids. This means that the velocity of the interfcae
is given by $(u\cdot \nu)\nu$. Since the normal component of
the velocity is continuous across the surface there is no ambiguity
in writing $u\cdot \nu$. The tangential components of $u_\pm$ need
not be continuous across $\Sigma (t)$, and indeed there may be jumps
in them. This allows possible slipping: the upper and
lower fluids moving in different directions tangent to $\Sigma (t)$.
Since only the normal component of the velocity vanishes at the
fixed upper and lower boundaries, $\{x_3=l\}$ and $\{x_3=-m\}$, the
fluids may also slip along the fixed boundaries.

To complete the statement of the problem, we must specify initial
conditions. We give the initial interface $\Sigma(0)=\Sigma_0$,
which yields the open sets $\Omega_\pm(0)$ on which we specify the
initial data for the density and the velocity
\begin{eqnarray*}  % \label{0108}
(\rho_\pm,u_\pm)(0,\cdot):\Omega_\pm(0)\rightarrow
({\mathbb{R}}^+,{\mathbb{R}}^3).
\end{eqnarray*}

To simply the equations we introduce the indicator function $\chi$
and denote
\begin{eqnarray*}  % \label{0109}
         \begin{array}{ll}
           \rho=\rho_+\chi_{\Omega_+}+\rho_-\chi_{\Omega_-},  \
   u=u_+\chi_{\Omega_+}+u_-\chi_{\Omega_-},\
p=p_+\chi_{\Omega_+}+p_-\chi_{\Omega_-}.
 \end{array}
\end{eqnarray*}
Hence the equations (\ref{0102}) are replaced by
\begin{equation} \label{0111}
\left\{
                  \begin{array}{ll}
                    \partial_t\rho+ {\mathrm{div}}(\rho u)=0,  \\[1mm]
                    \rho(\partial_t u+u\cdot\nabla u)+\nabla p(\rho)=-g\rho e_3+\rho
                    \omega u_2e_1-\rho
                    \omega u_1e_2,
                  \end{array}
                \right.
\end{equation}for $t>0$ and $x\in \Omega/\Sigma(t)$.
It will be convenient in our subsequent analysis to rewrite the momentum
equations in (\ref{0111}) by using the enthalpy function
\begin{eqnarray*} % \label{0112}
h(z)=\displaystyle\int_1^z\frac{p'(r)}{r}dr.
\end{eqnarray*}
The properties of $p$ ensure that $h\in C^\infty(0,\infty)$. Thus, (\ref{0111}) can be rewritten as
\begin{eqnarray}\label{0113}
\left\{
\begin{array}{ll}
 \partial_t\rho+ {\mathrm{div}}u=0,  \\[1mm]
\partial_tu+u\cdot\nabla u+\nabla h(\rho)=-ge_3+
                    2\omega u_2e_1-2\omega u_1e_2.
\end{array}
\right.
\end{eqnarray}

In the subsequent analysis it will be convenient to use the notation
\begin{eqnarray*}
\llbracket f \rrbracket:=f_+|_{x_3=0}-f_-|_{x_3=0}.
\end{eqnarray*}
\subsection{Steady-state solution}
\par \quad
In order to produce RTI, we first seek a steady-state solution with
$u_\pm=0$ and the interface given by $\{x_3=0\}$ for all $t\geq 0$.
Then $\Omega_+\equiv\Omega_+(t):={\mathbb{R}}^2\times (0,l)$,
$\Omega_-\equiv\Omega_-(t)=:{\mathbb{R}}^2\times (-m,0)$ for all
$t\geq 0$, and the equations reduce to the ODE
\begin{eqnarray}\label{0114}
\frac{d(p_\pm(\rho_\pm))}{d x_3}=-g\rho_\pm \ {\rm{in}} \
\Omega_\pm,
\end{eqnarray}
subject to the jump condition
\begin{eqnarray}\label{0115}
p_+(\rho_+)=p_-(\rho_-) \ {\rm{on}}\ \{x_3=0\}.
\end{eqnarray}
Such a solution depends only on $x_3$, so we consolidate notation by
assuming that $\rho_\pm$ are the restrictions to $(0,l)$ and
$(-m,0)$
 of a single function $\rho_0=\rho_0(x_3)$ that is smooth on $(-m,0)$ and $(0,l)$
 with a jump discontinuity across $\{x_3=0\}$.

From the assumption of the pressure function $p$ and the definition
of the set $Z$, there exist two positive constants $l$ and $m$ and
a solution $\rho_0$ to (\ref{0114})--(\ref{0115}) such that
\begin{itemize}
  \item $\rho_0$ is bounded above and below by positive constants on
  $(-m,l)$, and  $\rho_0$ is smooth when restricted to ($-m,0$)
  or $(0,l)$.
  \item $\llbracket \rho_0 \rrbracket=\rho_0^+|_{x_3=0}-\rho_0^-|_{x_3=0}>0$.
\end{itemize}
Please refer to \cite[Section 1.2]{Y-I} for more details concerning
construction of such a solution. \emph{In this paper, we always assume that
$l$, $m$ and the solution $\rho_0$ satisfy the above properties}.

\subsection{Formulation in Lagrangian coordinates}
 Since the movement of
the free interface $\Sigma(t)$ and the subsequent change of the
domains $\Omega_\pm(t)$ in Eulerian coordinates result in
mathematical difficulties, we switch our analysis to Lagrangian
coordinates as usual, so that the interface and the domains stay
fixed in time. To this end, we define the fixed Lagrangian domains
$\Omega_+={\mathbb{R}}^{2}\times (0,l)$ and
$\Omega_-={\mathbb{R}}^{2}\times (-m,0)$. We assume that there
exists invertible mappings
\begin{eqnarray*} % \label{0201}
\eta_\pm^0:\Omega_\pm\rightarrow \Omega_\pm(0),
\end{eqnarray*}
so that $\Sigma_0=\eta_+^0(\{x_3=0\})$,
$\{x_3=l\}=\eta_+^0(\{x_3=l\})$, and
$\{x_3=-m\}=\eta_-^0(\{x_3=-m\})$. The first condition means that
$\Sigma_0$ is parameterized by either of the mappings $\eta_+^0$
restricted to $\{x_3=0\}$, and the latter two conditions mean that
$\eta_\pm^0$ map the fixed upper and lower boundaries into
themselves. Define the flow maps $\eta_\pm$ by the solution to
\begin{eqnarray*}  %%%\label{0202}
\left\{
            \begin{array}{l}
\partial_t \eta_\pm(t,x)=u_\pm(t,\eta_\pm(t,x)),
\\
\eta_\pm(0,x)=\eta_\pm^0(x).
                  \end{array}    \right.
\end{eqnarray*}

We think of the Eulerian coordinates as $(t,y)$ with $y=\eta(t,x)$,
whereas we think of Lagrangian coordinates as the fixed $(t,x)\in
{\mathbb{R}}^+\times \Omega$, this implies that
$\Omega_{\pm}(t)=\eta_{\pm}(t,\Omega_{\pm})$ and that
$\Sigma(t)=\eta_+(t,\{x_3=0\})$, i.e., that the Eulerian domains of
upper and lower fluids are the image of $\Omega_\pm$ under the
mapping $\eta_\pm$ and that the free interface is parameterized by
$\eta_+(t,\cdot)$ restricted to $ {\mathbb{R}}^2\times\{0\}$. In
order to switch back and forth from Lagrangian to Eulerian
coordinates we assume that $\eta_\pm(t,\cdot)$ is invertible. Since
the upper and lower fluids may slip across one another, we must
introduce the slip map $S_\pm:$ ${\mathbb{R}}^+\times
{\mathbb{R}}^2\rightarrow {\mathbb{R}}^2\times\{0\}\subset
{\mathbb{R}}^2\times (-m,l)$ defined by
\begin{equation}\label{0203}
S_-(t,x')=\eta_-^{-1}(t,\eta_+(t,x',0)),\ x'\in {\mathbb{R}}^2
\end{equation}
and $S_+(t,\cdot)=S_-^{-1}(t,\cdot)$. The slip map $S_-$ gives the
particle in the lower fluid that is in contact with the particle of
the upper fluid at $x=(x_1,x_2,0)$ on the contact surface at time
$t$.

Setting $\eta=\chi_+\eta_++\chi_-\eta_-$, we define the Lagrangian
unknowns
\begin{eqnarray*} % \label{0204}
(v,q)(t,x)=(u,\rho)(t,\eta(t,x)),\ (t,x)\in {\mathbb{R}}^+\times
(\Omega/\{x_3=0\}).
\end{eqnarray*}
Defining the matrix $A:=(A_{ij})_{3\times 3}$ via
$A^T=(D\eta)^{-1}:=(\partial_j \eta_i)^{-1}_{3\times 3}$, and the
identity matrix $I=(I_{ij})_{3\times 3}$,  then in Lagrangian
coordinates the evolution equations for $\eta$, $v$,  $q$ are,
writing $\partial_j=\partial/\partial_{x_j}$,
\begin{equation}\label{0205} \left\{
                              \begin{array}{ll}
\partial_t \eta=v,\\
\partial_t q+q {\mathrm{tr}}(ADv)=0,\\[1mm]
\partial_t v+A\nabla h(q)=-gA\nabla\eta_3+2\omega v_2A\nabla\eta_1-2\omega v_1A\nabla\eta_2.
\end{array}
                            \right.
\end{equation}

Since the boundary jump conditions in Eulerian coordinates are phrased
in terms of jumps across the interface, the slip map must be employed
in Lagrangian coordinates. The jump conditions in Lagrangian coordinates are
\begin{eqnarray*}  %\label{0206}
\left\{
\begin{array}{ll}
(v_+(t,x',0)-v_-(S_-(t,x',0)))\cdot n(t,x',0)=0,\\
p(q_+(t,x',0))=p(q_-(S_-(t,x',0))),
\end{array}
\right.
\end{eqnarray*}
where we have written $n=\nu\circ \eta$, i.e.,
\begin{eqnarray*}
n:=\frac{\partial_1\eta_+\times\partial_2\eta_+}{|\partial_1\eta_+\times\partial_2\eta_+|}
\end{eqnarray*}
for the normal to the surface $\Sigma(t)=\eta_+(t,\{x_3=0\})$. Note that we
could just as well have phrased the jump conditions in terms of the
slip map $S_+$ and defined the interface and its normal vector in
terms of $\eta_-$. Finally, we require
\begin{eqnarray*}   %\label{0209}
v_-(t,x',-m)\cdot e_3=v_+(t,x',l)\cdot e_3=0.
\end{eqnarray*}
Note that since $\partial_t \eta=v$,
\begin{eqnarray*}  %%\label{0210}
e_3\cdot \eta_+(t,x',l)=e_3\cdot \eta_+^0(x',l)+\int_0^t e_3\cdot
v_+(t,x',l)\mathrm{d}s=l,
\end{eqnarray*}
which implies that $\eta_+(t,x',l)\in \{x_3=l\}$ for all $t\geq 0$,
i.e., that the part of the upper fluid in contact with the fixed boundary
$\{x_3=l\}$ never flows down from the boundary. It may, however,
slip along the fixed boundary, since we do not require
$v_+(t,x',l)\cdot e_i=0$ for $i=1$, $2$. A similar result holds
for $\eta_-$ at the lower fixed boundary $\{x_3=-m\}$.

In the steady-state case, the flow map is the identity mapping,
$\eta={\rm{Id}}$, so that $v=u$ and $q=\rho$. This means that the
steady-state solution, $\rho_0$, constructed above in Eulerian coordinates
 is also a steady state in Lagrangian coordinates. Now we
want to linearize the equations around the steady-state solution
$v=0$, $\eta={\rm{Id}}$, $q=\rho_0$, for which
$S_-={\rm{Id}}_{\{x_3=0\}}$ and $A=I$. The resulting linearized
equations read as
\begin{eqnarray}\label{1.2.14}
\left\{
\begin{array}{ll}
\partial_t\eta=v,\\
 \partial_tq+ \rho_0{\mathrm{div}}v=0,  \\[1mm]
                    \rho_0\partial_t v+\nabla (p'(\rho_0)q)=-gqe_3-g\rho_0\nabla\eta_3+2\rho_0
                    \omega v_2e_1-2\rho_0
                    \omega v_1e_2 .
\end{array}
\right.
\end{eqnarray}
The corresponding linearized jump conditions are
\begin{eqnarray*}  %\label{1.2.15}
\llbracket  v\cdot e_3 \rrbracket =0 \ {\mathrm{and}} \ \llbracket
p'(\rho_0)q\rrbracket=0,
\end{eqnarray*}
while the boundary conditions are
\begin{eqnarray*}  %\label{1.2.16}
v_-(t,x',-m)\cdot e_3=v_+(t,x',l)\cdot e_3=0.
\end{eqnarray*}
\subsection{Main results}
\par \quad
Before stating the first result concerning linear problem, we define some terms that will be
used throughout the paper.
 For a function $f\in L^2(\Omega)$, we define the horizontal
Fourier transform via
\begin{eqnarray*}  %\label{1.4.1}
\hat{f}(\xi_1, \xi_2, x_3)=\displaystyle\int_{{\mathbb{R}}^2}f(x_1,
x_2, x_3)e^{-i(x_1\xi_1+x_2\xi_2)}dx_1dx_2.
\end{eqnarray*}
By the Fubini and Parseval theorems, we have that
\begin{eqnarray*}  %\label{1.4.2}
\displaystyle\int_{\Omega}|f(x)|^2dx=\frac{1}{4\pi^2}\displaystyle\int_{\Omega}|\hat{f}(\xi,
x_3)|^2d\xi dx_3.
\end{eqnarray*}

For a function $f$ defined on $\Omega$ we write $f_+$ for the
restriction to $\Omega_+=\mathbb{R}^2\times (0,l)$ and $f_-$ for the
restriction to $\Omega_-=\mathbb{R}^2\times (-m,0)$. For $s\in
\mathbb{R}$, define the piecewise Sobolev space of order $s$ by
\begin{eqnarray}\label{1.4.3}
H^s(\Omega)=\left\{f|f_+\in H^s(\Omega_+), f_-\in
H^s(\Omega_-)\right\}
\end{eqnarray}
endowed with the norm
$\|f\|^2_{H^s}=\|f\|^2_{H^s(\Omega_+)}+\|f\|^2_{H^s(\Omega_-)}$. For
$k\in\mathbb{N}$ we can take the norms to be given by
\begin{equation*}  %\label{1.4.4}
\begin{aligned}
\|f\|^2_{H^k(\Omega_\pm)}:&=\sum_{j=0}^k\displaystyle\int_{\mathbb{R}^2\times
I_\pm}(1+|\xi|^2)^
{k-j}\left|\partial_{x_3}^j \hat{f}_\pm(\xi,x_3)\right|^2d\xi dx_3\\[2mm]
&=\sum_{j=0}^k\displaystyle\int_{\mathbb{R}^2}(1+|\xi|^2)^{k-j}\left|\partial_{x_3}^j
\hat{f}_\pm(\xi,\cdot)\right|^2_{L^2(I_\pm)}d\xi
\end{aligned}
\end{equation*}
for $I_-=(-m,0)$ and $I_+=(0,l)$. The main difference between the
piecewise Sobolev space $H^s(\Omega)$ and the ususal Sobolev space
is that we do not require functions in the piecewise Sobolev space
to have weak derivatives across the set $\{x_3=0\}$. Now, we may state our result on
the linear problem as follows:
\begin{Theorem}
The linear problem (\ref{1.2.14}) with the corresponding jump and
boundary conditions is ill-posed in the sense of Hadamard in
$H^k(\Omega)$ for every $k$. More precisly, for any
$k,j\in\mathbb{N}$ with $j\geq k$ and for any $T_0>0$ and $\alpha>0$,
there exists a sequence of solutions $\{(\eta_n,v_n,
q_n)\}^\infty_{n=1}$ to (\ref{1.2.14}), satisfying the corresponding
jump and boundary conditions, so that
\begin{eqnarray}\label{3.3.1}
\|\eta_n(0)\|_{H^j}+\|v_n(0)\|_{H^j}+\|q_n(0)\|_{H^j}\leq \frac1n ,
\end{eqnarray}
but
\begin{eqnarray}\label{3.3.2}
\|v_n(t)\|_{H^k}\geq \|\eta_n(t)\|_{H^k}\geq \alpha \ for \ all\ t\geq T_0.
\end{eqnarray}
\end{Theorem}

Theorem 1.1 shows discontinuous dependence on the initial data. More
precisely, there is a sequence of solutions with initial data
tending to $0$ in $H^k(\Omega)$, but which grow to be arbitrarily
large in $H^k(\Omega)$. Here we describe the framework
of the proof, which is inspired by \cite{Y-I}. First the resulting
linearized equations have coefficient functions that depend only on the
vertical variable, $x_3\in (-m,l)$. This allows us to seek ``normal
mode" solutions by taking the horizontal Fourier transform of the
equations and assuming that the solution grows exponentially in time by
the factor $e^{\lambda(|\xi|)t}$, where $\xi\in \mathbb{R}^2$ is the
horizontal spatial frequency and $\lambda(|\xi|)>0$. We show in
Theorem 2.2 that $\lambda(|\xi|)\rightarrow \infty$ in some unbounded
domain, the normal modes with a higher spatial frequency grow faster
in time, thus providing a mechanism for RTI.
 Indeed, we can form a Fourier synthesis of the normal
mode solutions constructed for each spatial frequency $\xi$ to
construct solutions of the linearized equations that grow
arbitrarily quickly in time, when measured in $H^k(\Omega)$ for any
$k\geq 0$. Comparing with \cite{Y-I}, here the main difficulty of
constructing this growing solutions lies in building the variational
structure of the linearized equations, because the natural
variational structure breaks down in the presence of the rotation
term. This difficulty will be circumvented in Section 2 by employing
an approach that was used first by Guo and Tice [4] to overcome a
similar difficulty arising from viscous compressible fluids, and
later adapted by Jiang, Jiang and Wang [7] in a nontrivial way to
construct growing mode solutions for viscous incompressible fluids
with magnetic field. Due to presence of the rotating term, we have
to impose stronger restrictions on the parameter $s$
and the spatial frequency $\xi$ in order to obtain the existence
result (\ref{2.2.12}) and the lower boundedness of $\lambda^2$ in (\ref{2.2.010})
for the rotating case. In addition, the auxiliary function
$\Phi(s)$ (see \cite[(3.52)]{Y-I2}), constructed by Guo and Tice to show
(\ref{2.2.12}), can not be applied to our case, and we have to construct a new
auxiliary function $F(s)$ in order to get (\ref{2.2.12}). At last, in
Section 3 we show a connection between the growth rate of solutions
to the linearized equations and the eigenvalues $\lambda(\xi)$,
which then gives rise to a uniqueness result (see Theorem 3.1). In
spite of the uniqueness, the linear problem is still ill-posed in
the sense of Hadamard in $H^k(\Omega)$ for any $k$, since solutions
do not depend continuously on the initial data.

With the linear ill-posedness established, we can obtain the ill-posedness of the
fully nonlinear problem in some sense. We rephrase the nonlinear equations (\ref{0205})
in a perturbation formulation around the steady state, that is,
$v=0,\ \eta=\eta^{-1}={\rm{Id}},\ q=\rho_0 $ with $A=I$ and $S_-=S_+={\rm{Id}}_{\{x_3=0\}}$.
Let
\begin{eqnarray*}  %\label{4.1.1}
\eta={\rm{Id}}+\tilde{\eta}, \ \eta^{-1}={\rm{Id}}-\zeta,\ v=0+v,\
q=\rho_0+\sigma,\  A=I-B,
\end{eqnarray*}
where
\begin{eqnarray*}  %\label{4.1.2}
B^T=\sum_{n=1}^{\infty}(-1)^{n-1}(D\tilde{\eta})^n.
\end{eqnarray*}

In order to deal with the term $h(q)=h(\rho_0+\sigma)$ we introduce the Taylor expansion
\begin{eqnarray*}  %\label{4.1.012}
h(\rho_0+\sigma)=h(\rho_0)+h'(\rho_0)\sigma+\mathfrak{R},
\end{eqnarray*}
where the remainder term is defined by
\begin{eqnarray*}  %\label{4.1.013}
\begin{array}{ll}
{\mathfrak{R}}(t,x)& =\displaystyle\int_0^{\sigma(t,x)}(\sigma(t,x)-z)h''(\rho_0(x)+z)dz\\[3mm]
& =\displaystyle\int_{\rho_0(x)}^{\rho_0(x)+\sigma(t,x)}
(\rho_0(x)+\sigma(t,x)-z)h''(z)dz .
\end{array}
\end{eqnarray*}
 Then (\ref{0205}) can be written for $\tilde{\eta},\ v,\ \sigma$ as
 \begin{eqnarray}\label{4.1.3}
 \left\{
 \begin{array}{ll}
 \partial_t\tilde{\eta}=v,\\[2mm]
 \partial_t\sigma+(\rho_0+\sigma)({\mathrm{div}} v-{\mathrm{tr}}(BDv))=0,\\[2mm]
\partial_t v+(I-B)\nabla(h'(\rho_0)\sigma+g\tilde{\eta}_3+{\mathfrak{R}})\\[2mm]
\qquad =2\omega v_2(I-B)\nabla({\rm{Id}}_1+\tilde{\eta}_1)-2\omega
v_1(I-B)\nabla({\rm{Id}}_2+\tilde{\eta}_2),
 \end{array}
 \right.
 \end{eqnarray}
 where $\mathrm{Id}_i$ denotes the $i$-th component of
 $\rm{Id}$, $i=1$, $2$.
We require the compatibility between $\zeta$ and $\tilde{\eta}$
given by
\begin{eqnarray}\label{4.1.4}
\zeta=\tilde{\eta}\circ({\rm{Id}}-\zeta).
\end{eqnarray}
The jump conditions across the interface are
\begin{eqnarray}\label{4.1.5}
\left\{
\begin{array}{ll}
(v_+(t,x',0)-v_-(t,S_-(t,x')))\cdot n(t,x',0)=0,\\[2mm]
p_+(\rho_0^++\sigma_+(t,x',0))=p_-(t,\rho_0^-+\sigma_-(S_-(t,x'))),
\end{array}
\right.
\end{eqnarray}
where the slip map (\ref{0203}) is rewritten as
\begin{eqnarray}\label{4.1.6}
S_-=({\rm{Id}}_{\mathbb{R}^2}-\zeta_-)\circ({\rm{Id}}_{\mathbb{R}^2}+\tilde{\eta}_+)
={\rm{Id}}_{\mathbb{R}^2}+\tilde{\eta}_+-\zeta_-\circ({\rm{Id}}_{\mathbb{R}^2}+\tilde{\eta}_+).
\end{eqnarray}

Finally, we require the boundary condition
\begin{eqnarray}\label{4.1.71}
v_-(t,x',-m)\cdot e_3=v_+(t,x',l)\cdot e_3=0.
\end{eqnarray}
We collectively refer to $(\ref{4.1.3})-(\ref{4.1.71})$ as the perturbed problem.
To shorten notation, for $k\geq 0$ we define
\begin{eqnarray*}  %\label{4.1.8}
\|(\tilde{\eta}, v, \sigma)(t)\|_{H^k}=\|\tilde{\eta}(t)\|_{H^k}+\|v(t)\|_{H^k}+\|\sigma(t)\|_{H^k}.
\end{eqnarray*}

In order to prove the ill-posedness for the perturbed problem by contradiction,
we state a definition introduced in \cite{Y-I}:
\begin{Definition} We say that the perturbed problem has property $EE(k)$ for
some $k\geq 3$ if there exist $\delta, t_0,C>0$ and
a function $F:[0,\delta)\rightarrow R^+$ satisfying
\begin{eqnarray*}  %\label{4.1.9}
\|(\tilde{\eta}_0, v_0,\sigma_0 )\|_{H^k}<\delta,
\end{eqnarray*}
there exist $(\tilde{\eta}, v, \sigma)\in
L^\infty((0,t_0),H^3(\Omega))$, such that
\begin{enumerate}
  \item[(1)] $(\tilde{\eta}, v, \sigma)(0)=(\tilde{\eta}_0, v_0, \sigma_0)$,
  \item[(2)] $\eta(t)={\rm{Id}}+\tilde{\eta}(t)$ is invertible and $\eta^{-1}(t)={\rm{Id}}-\zeta(t)$ for $0\leq t<t_0$,
  \item [(3)]$\tilde{\eta}$, $v$, $\sigma$, $\zeta$ solve the perturbed problem on $(0,t_0)\times\Omega$,
  \item[(4)]  we have the estimate
\begin{eqnarray*}  %\label{4.1.11}
\sup\limits_{0\leq t< t_0}\|(\tilde{\eta}, v, \sigma)(t)\|_{H^3}\leq
F(\|(\tilde{\eta}_0, v_0,\sigma_0)\|_{H^k}).
\end{eqnarray*}
\end{enumerate}
\end{Definition}
Here the $EE$ stands for existence and estimates, i.e. local-in-time
existence of solutions for small initial data, coupled to
$L^\infty(0,t_0;H^3(\Omega))$ estimates in terms of
$H^k(\Omega)$-norm of the initial data. If we were to add the
additional condition that such solutions be unique, then this trio
could be considered a well-posedness theory for the perturbed problem.

We can now show that property $EE(k)$ cannot hold for any $k\geq 3$. The proof utilizes
the Lipschitz structure of $F$ to show that property $EE(k)$
would give rise to certain estimates of solutions to the linearized
equations (\ref{1.2.14}) that cannot hold in general because of Theorem 1.1.

\begin{Theorem} The perturbed problem does not have property $EE(k)$ for any $k\geq 3$.
\end{Theorem}

\begin{Remark}
Theorem 1.1 and 1.2 show that rotating angular velocity $\boldsymbol{\omega}$ can not prevent
the linear and nonlinear RTI in the sense described in Theorem 1.1 and 1.2, respectively.
However, in the construction of the normal mode solution to the linearized system in section 2.2,
the rotation dose have a stabilizing effect on the growth rate $\lambda$ for sufficiently large fixed $|\xi|$.
 In fact, in section 2.2, we know that there exist a couple
 \begin{eqnarray*} (\varphi_0,\psi_0)\in {\mathcal{A}} =\left\{(\varphi,\psi)\in L^2(-m,l)
 \times H_0^1(-m,l)~\Big|~\frac12\displaystyle\int_{-m}^l\rho_0(\varphi^2+\psi^2)dx_3=1\right\}
\end{eqnarray*} such that  \begin{eqnarray*} &&0>-\lambda^2=
\frac{1}{2}\displaystyle\int_{-m}^l\left\{p'(\rho_0)\rho_0(\psi'_0+|\xi|
\varphi_0)^2-2g|\xi|\rho_0\psi_0\varphi_0+4\omega^2
\rho_0\varphi^2_0/\lambda^{2}\right\}dx_3\\
&&\qquad\qquad  = \inf\limits_{(\varphi,\psi)\in
\mathcal{A}}\frac{1}{2}\displaystyle\int_{-m}^l\left\{p'(\rho_0)\rho_0(\psi'+|\xi|
\varphi)^2-2g|\xi|\rho_0\psi\varphi+4\omega^2
\rho_0\varphi^2_0/\lambda^{2}\right\}dx_3.
 \end{eqnarray*}Obviously, $\varphi_0\equiv\!\!\!\!\!/\ 0$.
On the other hand, we denote the growth rate by $\lambda_0$ for the
equations of compressible inviscid fluids without  rotation (i.e.,
$\boldsymbol{\omega}=\textbf{0}$), that is
 \begin{eqnarray*}
&&0>-\lambda^2_0= \inf\limits_{(\varphi,\psi)\in
\mathcal{A}}\frac{1}{2}\displaystyle\int_{-m}^l\left\{p'(\rho_0)\rho_0(\psi'+|\xi|
\varphi)^2-2g|\xi|\rho_0\psi\varphi\right\}dx_3.
 \end{eqnarray*}
 Thus we have
 \begin{eqnarray*}
&&-\lambda^2\geq -\lambda^2_0+2\int_{-m}^l\omega^2
\rho_0\varphi^2_0/\lambda^{2}dx_3>-\lambda^2_0\mbox{ for }\omega\neq
0, \end{eqnarray*}which implies $\lambda<\lambda_0$.
\end{Remark}

\section{Construction of a growing solution to (\ref{1.2.14})}
\setcounter{equation}{0}

\quad \ We wish to construct a solution to the linearized equations (\ref{1.2.14}) that has a growing $H^k$-norm
for any $k$. We will construct such solutions via Fourier synthesis
by first constructing a growing mode for a fixed spacial frequency.
\subsection{Growing mode and Fourier transform}
\quad To begin, we make an ansatz
\begin{eqnarray*}  %\label{2.1.1}
v(t,x)=w(x) e^{\lambda t},\ \ q(t,x)=\tilde{q}(x)e^{\lambda t}, \ \eta(t,x)=\tilde{\eta}(x)e^{\lambda t},
\end{eqnarray*}
for some $\lambda>0$. Substituting this ansatz into (\ref{1.2.14}), eliminating $\tilde{\eta}$ and $ \tilde{q}$ by using the first two equations,
we arrive at the
time-invariant system for $w=(w_1,w_2,w_3)$:
\begin{eqnarray}\label{2.1.2}
\lambda^2\rho_0w-\nabla(p'(\rho_0)\rho_0{\mathrm{div}}w)=g\rho_0{\mathrm{div}}w
e_3-g\rho_0\nabla w_3+2\rho_0\omega w_2 e_1-2\rho_0 \omega w_1 e_2
\end{eqnarray}
with the corresponding jump conditions
\begin{eqnarray*}  %\label{2.1.3}
\llbracket   w_3\rrbracket=0 \ {\mbox{ and }}\  \llbracket
p'(\rho_0)\rho_0{\mathrm{div}} w\rrbracket=0,
\end{eqnarray*}
and the boundary conditions
\begin{eqnarray*}  %\label{2.1.4}
w_3(t,x',-m)=w_3(t,x',l)=0,\ x'\in \mathbb{R}^2.
\end{eqnarray*}

Since the jump only occurs in the $e_3$ direction, we are free to
take the horizontal Fourier transform, which we denote with either
$\hat{\cdot}$ or $\mathcal{F}$, to reduce to a system of ODEs in
$x_3$ for each fixed spacial frequency.

We take the horizontal Fourier transform of $w_1$, $w_2$, $w_3$ in
(\ref{2.1.2}) and fix a spacial frequency $\xi=(\xi_1, \xi_2)\in
\mathbb{R}^2$. Define the new unknowns
$$\varphi(x_3)=i\hat{w}_1(\xi_1,\xi_2,x_3),\;\; \theta(x_3)=i
\hat{w}_2(\xi_1,\xi_2,x_3),\;\; \psi(x_3)=\hat{w}_3(\xi_1,\xi_2,x_3). $$
To write down the equations for $\varphi$, $\theta$, $\psi$, we denote
$'=d/dx_3$ to arrive at the following system of ODEs
\begin{eqnarray}\label{2.2.1}
\left\{
\begin{array}{ll}
\lambda^2\rho_0\varphi+\xi_1[p'(\rho_0)\rho_0(\xi_1\varphi+\xi_2\theta+\psi')]=\xi_1g\rho_0\psi+2\rho_0\omega\theta,\\[2mm]
\lambda^2\rho_0\theta+\xi_2[p'(\rho_0)\rho_0(\xi_1\varphi+\xi_2\theta+\psi')]=\xi_2g\rho_0\psi-2\rho_0\omega\varphi,\\[2mm]
\lambda^2\rho_0\psi-[p'(\rho_0)\rho_0(\xi_1\varphi+\xi_2\theta+\psi')]'=g\rho_0(\xi_1\varphi+\xi_2\theta),
\end{array}
\right.
\end{eqnarray}
 along with the jump conditions
 \begin{eqnarray*}  %\label{2.2.2}
 \left\{
\begin{array}{ll}
\llbracket \psi \rrbracket=0,\\[2mm]
\llbracket
p'(\rho_0)\rho_0(\xi_1\varphi+\xi_2\theta+\psi')\rrbracket=0,
\end{array}
\right.
\end{eqnarray*}
and boundary conditions
\begin{eqnarray}\label{2.2.3}
\psi(-m)=\psi(l)=0.
\end{eqnarray}

We can reduce the complexity of the problem by removing the
component $\theta$. For that purpose, note that if $\varphi$, $\theta$,
$\psi$ solve the above equations for $\xi_1$, $\xi_2$ and $\lambda$,
then for any rotation operator $\mathcal{R}\in {\mathcal{SO}(2)}$,
$(\tilde{\varphi},\tilde{\theta}):=\mathcal{R}(\varphi,\theta)$
solves the same equations for
$(\tilde{\xi}_1,\tilde{\xi}_2):={\mathcal{R}}(\xi_1,\xi_2)$ with
$\psi$, $\lambda$ unchanged. So, by choosing an appropriate
rotation, we may assume without loss of generality that $\xi_2=0$
and $\xi_1=|\xi|\geq 0$. In this setting, $\theta$ solves
\begin{eqnarray*}  %\label{2.2.03}
\lambda^2\rho_0\theta=-2\rho_0\omega\varphi.
\end{eqnarray*}
Putting this identity into (\ref{2.2.1}), we arrive at
\begin{eqnarray}\label{2.2.04}
\left\{
\begin{array}{ll}
\lambda^2\rho_0\varphi+|\xi|[p'(\rho_0)\rho_0(|\xi|\varphi+\psi')]
=|\xi|g\rho_0\psi-4\omega^2{\lambda^{-2}}{\rho_0}\varphi,
\\[2mm]
\lambda^2\rho_0\psi-[p'(\rho_0)\rho_0(|\xi|\varphi+\psi')]'=|\xi|g\rho_0\varphi,
\end{array}
\right.
\end{eqnarray}
along with the jump conditions
\begin{eqnarray}\label{2.2.040}
 \left\{
\begin{array}{ll}
\llbracket \psi \rrbracket=0,\\[2mm]
\llbracket p'(\rho_0)\rho_0(|\xi|\varphi+\psi') \rrbracket=0,
\end{array}
\right.
\end{eqnarray}
and boundary conditions (\ref{2.2.3}).

In the absence of rotation ($\omega=0$) and for a fixed spatial
frequency $\xi\neq 0$, the equations (\ref{2.2.04}), (\ref{2.2.040})
and (\ref{2.2.3}) can be viewed as an eigenvalue problem with
eigenvalue $-\lambda^2$. Such a problem has a natural variational
structure that allows us to construct solutions via the direct
method and a variational characterization of the eigenvalues via
\begin{eqnarray*} %\label{2.2.041}
-\lambda^2=\inf\frac{E(\varphi,\psi)}{J(\varphi,\psi)},
\end{eqnarray*}
where
\begin{eqnarray*} %\label{2.2.042}
E(\varphi,\psi)=\frac{1}{2}\displaystyle\int_{-m}^l
p'(\rho_0)\rho_0(\psi'+|\xi| \varphi)^2-2g|\xi|\rho_0\psi\varphi
dx_3
\end{eqnarray*}
and
\begin{eqnarray*}  %\label{2.2.043}
J(\varphi,\psi)=\frac12\displaystyle\int_{-m}^l\rho_0(\varphi^2+\psi^2)dx_3.
\end{eqnarray*}
This variational structure was essential to the analysis in
\cite{Y-I}, where the ill-posedness results for both the inviscid
linearized problem and the inviscid non-linear problem ((\ref{4.1.3})
with $\omega=0$) were shown. Unfortunately, when rotation is
present, the natural variational structure breaks down. In order to
circumvent this problem and restore the ability to use the variational
method, first we artificially remove the dependence of (\ref{2.2.04}) on ${\lambda^{-2}}$ by
 defining $s:={\lambda^{-2}}>0$, and then consider a family
($s>0$) of the modified problems given by

\begin{eqnarray}\label{2.2.4}
\left\{
\begin{array}{ll}
-\lambda^2\rho_0\varphi=|\xi|[p'(\rho_0)\rho_0(|\xi|\varphi+\psi')]-|\xi|g\rho_0\psi+4\omega^2s\rho_0 \varphi,\\[2mm]
-\lambda^2\rho_0\psi=-[p'(\rho_0)\rho_0(|\xi|\varphi+\psi')]'-|\xi|g\rho_0\varphi,
\end{array}
\right.
\end{eqnarray}
along with the jump conditions
\begin{eqnarray}\label{2.2.5}
 \left\{
\begin{array}{ll}
\llbracket \psi \rrbracket=0,\\[2mm]
\llbracket p'(\rho_0)\rho_0(|\xi|\varphi+\psi')\rrbracket=0,
\end{array}
\right.
\end{eqnarray}
and boundary conditions (\ref{2.2.3}).

\subsection{Construction of a solution to (\ref{2.2.4})}
\par \quad
In the following, we establish the variation framework by formulating constrained minimization.
Multiplying $(\ref{2.2.4})_1$  and $(\ref{2.2.4})_2$ by $\varphi$ and $\psi$, respectively,
 we add the resulting equations, integrate over
  $(-m,l)$, integrate by parts, and apply the boundary and jump conditions to deduce that
\begin{eqnarray*}  %\label{2.3.1}
-\frac{\lambda^2}{2}\displaystyle\int_{-1}^1\rho_0(\varphi^2+\psi^2)dx_3
=\frac{1}{2}\displaystyle\int_{-1}^1p'(\rho_0)\rho_0(\psi'+|\xi|
\varphi)^2-2g|\xi|\rho_0\psi\varphi+4\omega^2 s\rho_0\varphi^2dx_3.
\end{eqnarray*}
Notice that for any fixed $\xi$, this is a standard eigenvalue problem for $-\lambda^2$.
 It allows us to use the variational method to construct solutions.
To this end, we define the energies
\begin{eqnarray*}  %\label{2.3.2}
E(\varphi,\psi)=\frac{1}{2}\displaystyle\int_{-m}^l\left\{p'(\rho_0)\rho_0(\psi'+|\xi|
\varphi)^2-2g|\xi|\rho_0\psi\varphi+4\omega^2
s\rho_0\varphi^2\right\}dx_3
\end{eqnarray*}
and
\begin{eqnarray*}  %\label{2.3.3}
J(\varphi,\psi)=\frac12\displaystyle\int_{-1}^1\rho_0(\varphi^2+\psi^2)dx_3,
\end{eqnarray*}
which are both well-defined on the space $L^2(-m,l)\times
H_0^1(-m,l)$.
 We define the set
\begin{eqnarray*}  %\label{2.3.4}
{\mathcal{A}} =\{(\varphi,\psi)\in L^2(-m,l)\times
H_0^1(-m,l)~|~J(\varphi,\psi)=1\},
\end{eqnarray*}

We want to show that the infimum of $E(\varphi,\psi)$ over the set
${\mathcal{A}}$ can be achieved and is negative, and that the minimizer
solves the problem $(\ref{2.2.4}), (\ref{2.2.5})$ and
(\ref{2.2.3}). Notice that by employing the identity
$-2ab=(a-b)^2-(a^2+b^2)$ and the constraint on $J(\varphi,\psi)$, we may rewrite
\begin{equation}\label{2.3.04}
\begin{aligned}E(\varphi,\psi)&=-g|\xi|+\frac12\displaystyle\int_{-m}^lp'(\rho_0)\rho_0(\psi'+|\xi|\varphi)^2
+g|\xi|\rho_0(\varphi-\psi)^2dx_3
\\[2mm]
&\quad +2\omega^2{s}\displaystyle\int_{-m}^l\rho_0\psi^2dx_3\geq -g|\xi| ,
\qquad (\varphi,\psi)\in {\mathcal{A}}.
\end{aligned}
\end{equation}
In order to emphasize the dependence on $s\in(0,\infty)$, we will sometimes write
\begin{eqnarray*}  %\label{2.3.041}
E(\varphi,\psi)=E(\varphi,\psi;s)
\end{eqnarray*}
and
 \begin{eqnarray}\label{2.3.5111111}
\mu(s)=:-\lambda^2(s):= \inf\limits_{\psi\in \mathcal{A}}
E(\varphi,\psi;s).
 \end{eqnarray}
By inequality (\ref{2.3.04}) we have
\begin{eqnarray}\label{2.3.512223}\lambda\leq \sqrt{g|\xi|}. \end{eqnarray}
%%%%
\begin{Proposition}\label{pro:0201}
There exist three constants $R_1\geq 2$,  $C_0$ and $C_1$ depending
on the quantities $\rho_0^\pm$, $p_\pm$, $g$, $l$, $m$, $\omega$, so
that
  \begin{eqnarray}\label{2.2.10}
\mu(s)\leq -C_0|\xi|+sC_1 \mbox{ for any }|\xi|\geq R_1.
\end{eqnarray}
In particular, we have
  \begin{eqnarray*}
\mu(s)< 0\mbox{ for any  }|\xi|\geq R_1\mbox{ and }s<C_0|\xi|/C_1.
\end{eqnarray*}
\end{Proposition}
{\bf{Proof}}. Since both $E$ and $J$ are homogeneous of degree 2, it
suffices to show that
\begin{eqnarray*}  %\label{2.2.6}
\inf\limits_{(\varphi,\psi)\in L^2\times
H_0^1}\frac{E(\varphi,\psi)}{J(\varphi,\psi)} <0.
\end{eqnarray*}
 We assume that $\varphi=-\psi'/|\xi|$, such
that the first integrand term in $E(\varphi,\psi)$ vanishes, that is,
\begin{equation}\label{2.2.7}
\begin{aligned}
E\left(-\frac{\psi'}{|\xi|},\psi\right)&=\displaystyle\int_{-m}^lg\rho_0\psi\psi'dx_3+
\frac{2\omega^2 s}{|\xi|^2}\displaystyle\int_{-m}^l\rho_0\psi'^2dx_3\\[2mm]
&=-\frac{g\llbracket \rho_0 \rrbracket\psi^2(0)}{2} -\frac{1}{2}
\displaystyle\int_{-m}^l g\rho_0'\psi^2dx_3
+\frac{2\omega^2 s}{|\xi|^2}\displaystyle\int_{-m}^l\rho_0\psi'^2dx_3\\[2mm]
&=-\frac{g\llbracket \rho_0 \rrbracket\psi^2(0)}{2}+\frac{g^2}{2}
\displaystyle\int_{-m}^l\frac{\rho_0}{p'(\rho_0)}\psi^2dx_3
+\frac{2\omega^2s}{|\xi|^2}\displaystyle\int_{-m}^l\rho_0\psi'^2dx_3 ,
\end{aligned}
\end{equation}
where we have used the fact that $\rho_0$ solves (\ref{0114}).
Notice that $\llbracket \rho_0 \rrbracket>0$ so that the right-hand side is not positive definite.
From (\ref{2.2.7}) we see that the rotation term diminishes obviously the growth of instability.

For $|\xi|\geq 2$ we define the test function $\psi_{|\xi|}\in
H_0^1(-m,l)$ defined by
\begin{eqnarray}\label{2.2.8}
\psi_{|\xi|}(x_3)=\left\{
\begin{array}{ll}
\left(1-\frac{x_3}{l}\right)^{\frac{|\xi|}{2}},\ x_3\in[0,l)\\[2mm]
\left(1+\frac{x_3}{m}\right)^{\frac{|\xi|}{2}},\ x_3\in(-m,0).
\end{array}
\right.
\end{eqnarray}
Thus, we can estimate that
\begin{eqnarray}\label{2.2.9}
E\left(-\frac{\psi'}{|\xi|},\psi\right)\leq-\frac{g\llbracket \rho_0 \rrbracket}{2}
+\frac{A_1(1+{s})}{1+|\xi|}
\end{eqnarray}
and
\begin{eqnarray}\label{2.2.0110123}
\frac{A_2}{1+|\xi|} \leq J\left(-\frac{\psi'_{|\xi|}}{|\xi|}, \psi_{|\xi|}\right)\leq \frac{A_3
}{1+|\xi|}
\end{eqnarray}
for some constants $A_1$, $A_2$ and $A_3>0$ depending on $\rho_0^\pm$, $p_\pm$,
$l$, $m$, $g$ and $\omega$, but not on $|\xi|$.
Making use of (\ref{2.2.9}) and (\ref{2.2.0110123}), we infer that
\begin{eqnarray}\label{112.2.0111}
\frac{E(-\psi'_{|\xi|}/|\xi|, \psi_{|\xi|})}{J(-\psi'_{|\xi|}/|\xi|, \psi_{|\xi|})}\leq A_4-A_5|\xi|+
C_1 s
\end{eqnarray}
for some positive constants  $C_1$, $A_4$ and $A_5$ depending on
the same parameters, but not on $|\xi|$.

Obviously, there exist a sufficiently large positive
constant $R_1>0$ and a positive constant $C_0$ depending on $A_4$ and $A_5$, such that
\begin{equation}\label{11212.2.0111}A_4-A_5|\xi|/2\leq 0\mbox{ for any } |\xi|\geq R_1.
\end{equation}
Inserting (\ref{11212.2.0111}) into (\ref{112.2.0111}), we obtain
\begin{eqnarray*}
\frac{E(-\psi'_{|\xi|}/|\xi|, \psi_{|\xi|})}{J(-\psi'_{|\xi|}/|\xi|, \psi_{|\xi|})}\leq -C_0 |\xi|+
C_1 s\mbox{ for }C_0=A_5/2,
\end{eqnarray*} which implies the desired conclusion.
\hfill $\Box$
\vspace{2mm}

The key point in the proof argument of Proposition \ref{pro:0201} is to construct
a pair $(\varphi,\psi)$, such that
\begin{eqnarray*}
-\frac{g\llbracket \rho_0 \rrbracket\psi^2(0)}{2}<0,
\end{eqnarray*}
which in particular requires that $\psi(0)\neq 0$. We can show that this property
is satisfied by any $(\varphi,\psi)\in{\mathcal{A}}$ with $E(\varphi,\psi)<0$.
%%%%%%%%%%%%%%%%%%%%%%%%%%%%%%%%%%%%%%
\begin{Lemma}\label{lem:0201}
Suppose that $(\varphi,\psi)\in {\mathcal{A}}$ satisfies $E(\varphi, \psi)<0$. Then $\psi(0)\neq 0$.
\end{Lemma}
{\bf{Proof}}. A completion of square allows us to write
$$ p'(\rho_0)\rho_0(\psi'+|\xi|\varphi)^2-2g\rho_0|\xi|\psi
=\left(\sqrt{p'(\rho)\rho_0}(\psi'+|\xi|\varphi)
-\frac{g\sqrt{\rho_0}\psi}{\sqrt{p'(\rho_0)}}\right)^2+2g\rho_0\psi\psi'-\frac{g^2\rho_0}{p'(\rho_0)}\psi^2.
$$
Thus, similarly to (\ref{2.2.7}), we can rewrite the energy as
\begin{eqnarray*}
E(\varphi,\psi)=-\frac{g\llbracket \rho_0
\rrbracket\psi^2(0)}{2}+\frac12\displaystyle\int_{-m}^l\left(\sqrt{p'(\rho)
\rho_0}(\psi'+|\xi|\varphi)-\frac{g\sqrt{\rho_0}}{\sqrt{p'(\rho_0)}}\psi\right)^2dx_3
+2{\omega^2s}\displaystyle\int_{-m}^l\rho_0\varphi^2dx_3,
\end{eqnarray*}
from which we deduce that if $E(\varphi,\psi)<0$, then $\psi(0)\neq 0$. \hfill $\Box$
\vspace{2mm}

From now on, we denote
\begin{eqnarray*}
\mathcal{M}:=\{(\xi,s)~|~\xi\in \mathbb{R}^2,\ |\xi|\geq R_1,\ s<C_0|\xi|/C_1\},
\end{eqnarray*}where the constants $R_1$, $C_0$ and $C_1$ are  from Proposition \ref{pro:0201}.
Next, we can show that a minimizer exists and that the minimizer satisfies
(\ref{2.2.4}), (\ref{2.2.5}) and (\ref{2.2.3}) for each
$(\xi, s)\in \mathcal{M}$ in the same manner as in \cite{Y-I2}.

\begin{Proposition}
For any fixed $(\xi, s)\in \mathcal{M}$,  $E$ achieves its infinimum on $\mathcal{A}$.
\end{Proposition}
{\bf{Proof}}. (\ref{2.3.04}) shows that $E$ is bounded from below on
${\mathcal{A}}$. Let $(\varphi_n,\psi_n)\in {\mathcal{A}}$ be a
minimizing sequence. Then $\varphi_n$ is bounded in $L^2(-m,l)$ and $\psi_n$
 is bounded in $H_0^1(-m,l)$, so up to the extraction of a subsequence
$\varphi_n\rightharpoonup \varphi$ weakly in $L^2$, $\psi_n\rightharpoonup \psi$ weakly in $H_0^1$,
and $\psi_n\rightarrow \psi$ strongly in $L^2$. In view of
the weak lower semi-continuity and the strong $L^2$-convergence $\psi_n\rightarrow\psi$, we find that
\begin{eqnarray*}  %\label{2.3.05}
E(\varphi,\psi)\leq \lim\inf\limits_{n\rightarrow \infty}E(\varphi_n,\psi_n)=\inf\limits_{\mathcal{A}} E.
\end{eqnarray*}
All that remains is to show $(\varphi,\psi)\in{\mathcal{A}}$.

Again by the lower semi-continuity, we see that $J(\varphi,\psi)\leq 1$.
Suppose by contradiction that $J(\varphi,\psi)<1$. By the
homogeneity of $J$ we may find $\alpha>1$ so that $J(\alpha\varphi,\alpha \psi)=1$,
i.e., we may scale up $(\varphi,\psi)$ so that
$(\alpha\varphi,\alpha \psi)\in{\mathcal{A}}$. By Proposition 2.1
we know that $\inf\limits_{\mathcal{A}}E<0$, from which we deduce that
\begin{eqnarray*}
E(\alpha\varphi,\alpha \psi)=\alpha^2E(\varphi, \psi)\leq \alpha^2
\inf\limits_{\mathcal{A}}E<\inf\limits_{\mathcal{A}}E,
\end{eqnarray*}
which is a contradiction since $(\alpha\varphi,\alpha \psi)\in{\mathcal{A}}$.
Hence $J(\varphi,\psi)=1$ so that $(\varphi,\psi)\in{ \mathcal{A}}$.
\hfill $\Box$

\begin{Proposition}\label{pro:0203} Let $(\xi, s)\in \mathcal{M}$,
and  $(\varphi,\psi)\in\mathcal{A}$ be the minimizer of $E$ constructed in Proposition 2.2.
Let $\mu:=-\lambda^2:=E(\varphi,\psi)$. Then $\varphi,\psi$ satisfy
\begin{eqnarray}\label{2.3.051}
\left\{
\begin{array}{ll}
\mu\rho_0\varphi=|\xi|[p'(\rho_0)\rho_0(|\xi|\varphi+\psi')]-|\xi|g\rho_0\psi+4\omega^2 s\rho_0\varphi,\\[2mm]
\mu\rho_0\psi=-[p'(\rho_0)\rho_0(|\xi|\varphi+\psi')]'-|\xi|g\rho_0\varphi
\end{array}
\right.
\end{eqnarray}
along with the jump conditions (\ref{2.2.040}) and boundary
conditions (\ref{2.2.3}). Moreover, the solutions are smooth when
restricted to either $(-m,0)$ or $(0,l)$.
\end{Proposition}
{\bf{Proof}}. Fix $(\varphi_0,\psi_0)\in L^2(-m,l)\times
H_0^1(-m,l)$. Define
\begin{eqnarray*}  %\label{2.3.052}
j(t,r)=J(\varphi+t\varphi_0+r\varphi,\psi+t\psi_0+r\psi )
\end{eqnarray*}
and note that $j(0,0)=1$. Moreover, $j$ is smooth, and
$$  %% \label{2.3.053}
\frac{\partial j}{\partial t}(0,0)=\displaystyle\int_{-m}^l\rho_0[\varphi_0\varphi+\psi_0\psi]dx_3,
\qquad %% \label{2.3.054}
\frac{\partial j}{\partial r}(0,0)=\displaystyle\int_{-m}^l\rho_0(\varphi^2+\psi^2)=2.  $$
Thus, by the inverse function theorem, we can find a $C^1$-function $r=\sigma(t)$
in a neighborhood of $0$, such that $\sigma(0)=0$ and $j(t,\sigma(t))=1$.
We may differentiate the last equation to infer that
\begin{eqnarray*}
\frac{\partial j}{\partial t}(0,0)+\frac{\partial j}{\partial l}(0,0)\sigma'(0)=0,
\end{eqnarray*}
whence,
\begin{eqnarray*}
\sigma'(0)=-\frac12\frac{\partial j}{\partial
t}(0,0)=-\frac12\displaystyle\int_{-m}^l\rho_0[\varphi_0\varphi+\psi_0\psi]dx_3.
\end{eqnarray*}

Since $(\varphi, \psi)$ is a minimizer of $E$ over ${\mathcal{A}}$, one has
\begin{eqnarray*}
0=\frac{d}{dt}\bigg|_{t=0}E(\varphi+t\varphi_0+\sigma(t)\varphi,\psi+t\psi_0+\sigma(t)\psi
),
\end{eqnarray*}
 which implies that
 \begin{eqnarray*}
  0 &=& \displaystyle\int_{-m}^lp'(\rho_0)\rho_0(\psi'+|\xi|\varphi)
 (\psi_0'+\sigma'(0)\psi'+|\xi|\varphi_0+|\xi|\sigma'(0)\varphi) \\
 && -g|\xi|\rho_0(\psi(\varphi_0+\sigma'(0)\varphi)
 +\varphi(\psi_0+\sigma'(0)\psi))+4\omega^2s\rho_0 \varphi(\varphi_0+\sigma'(0)\varphi)dx_3.
  \end{eqnarray*}
The above equation, by rearranging and plugging in the value of $\sigma'(0)$, can be rewrite as
\begin{equation*}
\begin{aligned}
&\int_{-m}^l
p'(\rho_0)\rho_0(\psi'+|\xi|\varphi)(\psi_0'+|\xi|\varphi_0)-g|\xi|
\rho_0(\psi\varphi_0+\varphi\psi_0)+4\omega^2s\rho_0\varphi\varphi_0dx_3
\\
&=\mu\displaystyle\int_{-m}^l\rho_0[\varphi_0\varphi+\psi_0\psi]dx_3,
\end{aligned}\end{equation*}
where the lagrange multiplier (eigenvalue) is $\mu=E(\varphi,\psi)$.

By making variations with $\varphi_0,\psi_0$ compactly supported in
either $(-m,0)$ or $(0,l)$, we find that $\varphi$ and $\psi$
satisfy the equations (\ref{2.3.051}) in the weak sense in $(-m,0)$
and $(0,l)$. Standard bootstrapping arguments then show that
$(\varphi,\psi)$ are in $H^k(-m,0)$ (resp. $H^k(0,l)$) for all
$k\geq 0$ when restricted to $(-m,0)$ (resp. $(0,l)$), and hence the
functions are smooth when restricted to either interval. This
implies that the equations are also classically satisfied on
$(-m,0)$ and $(0,l)$. Since $(\varphi,\psi)\in H^2$, the traces of
the functions and their derivatives are well-defined at the
endpoints $x_3=-m$, $0$ and $l$, and it remains to show that the
jump conditions are satisfied at $x_3=0$ and the boundary conditions
satisfied at $x_3=-m$ and $l$. Making variations with respect to
arbitrary $\varphi_0, \psi_0\in C_c^\infty(-m,l)$, we find that the
jump condition
\begin{eqnarray*}
\llbracket   p'(\rho_0)\rho_0(|\xi|\varphi+\psi')\rrbracket=0.
\end{eqnarray*}
must be satisfied. Note that the conditions
 $\llbracket \psi \rrbracket =0$ and $\psi(-m)=\psi(l)=0$ are satisfied
 trivially since $\psi\in H_0^1(-m,l)\hookrightarrow C_0^{0,1/2}(-m,l)$. \hfill $\Box$
 \vspace{2mm}

The next result establishes continuity property of the eigenvalue $\mu(s)$
which will be used in finding a $s$ with $s=1/\lambda^2(s)$ in Theorem 2.1.
%%%%%%%%%%%%%%%%%%%%%%%
\begin{Proposition}Let $\mu :(0,\infty)\rightarrow {\mathbb{R}}$ be given by (\ref{2.3.5111111}),
then $\mu\in C_{{\rm{loc}}}^{0,1}(0,\infty)$, and in particular, $\mu\in C^0(0,\infty)$.
\end{Proposition}
{\bf{Proof.}} Fix a compact interval $I=[a,b]\subset(0,\infty)$, and
fix any pair $(\varphi_0,\psi_0)\in \mathcal{A}$. We may decompose
$E$ according to
\begin{eqnarray}\label{2.2.0112}
E(\varphi,\psi;s)=E_0(\varphi,\psi)+sE_1(\varphi,\psi),
\end{eqnarray}
where
\begin{eqnarray*}   &&
E_0(\varphi,\psi):=\frac{1}{2}\displaystyle\int_{-m}^l\left\{p'(\rho_0)\rho_0(\psi'+|\xi|
\varphi)^2-2g|\xi|\rho_0\psi\varphi\right\}dx_3, \\
&& E_1(\varphi,\psi):=2{\omega^2}\displaystyle\int_{-m}^l\rho_0\psi^2 dx_3\geq 0.
\end{eqnarray*}
The non-negativity of $E_1$ implies that $E$ is non-decreasing in $s$ with fixed $(\varphi,\psi)\in {\mathcal{A}}$.

Now, by Proposition 2.2, for each $s\in(0,\infty)$ we can find a pair
$(\varphi_s,\psi_s)\in {\mathcal{A}}$ so that
\begin{eqnarray*}  %\label{2.2.0115}
E(\varphi_s,\psi_s;s)=\inf\limits_{(\varphi,\psi)\in{\mathcal{A}}}E(\varphi,\psi;s)=\mu(s).
\end{eqnarray*}
We deduce from the non-negativity of $E_1$, the minimality of $(\varphi_s, \psi_s)$
and the equality (\ref{2.3.04}) that
\begin{eqnarray*}  %\label{2.2.0116}
E(\varphi_0,\psi_0;b)\geq E(\varphi_0,\psi_0;s)\geq E(\varphi_s,\psi_s;s)\geq s E_1(\varphi_s,\psi_s)-g|\xi|
\end{eqnarray*}
for all $s\in Q$, whichi implies that there exists a constant
$0<K=K(a,b,\varphi_0,\psi_0, g, |\xi|)<\infty$, such that
\begin{eqnarray}\label{2.2.0117}
\sup\limits_{s\in Q}E_1(\varphi_s,\psi_s)\leq K.
\end{eqnarray}

Let $s_i\in Q$ for $i=1,2.$ Using the minimality of $(\varphi_{s_1},\psi_{s_1})$
compared to $(\varphi_{s_2},\psi_{s_2})$, we see that
\begin{eqnarray*}  %\label{2.2.0118}
\mu(s_1)=E(\varphi_{s_1},\psi_{s_1};s_1)\leq E(\varphi_{s_2},\psi_{s_2};s_1).
\end{eqnarray*}
Recalling the decomposition (\ref{2.2.0112}), we can bound
\begin{eqnarray*} %\label{2.2.0119}
\begin{array}{ll}
E(\varphi_{s_2},\psi_{s_2};s_1)&\leq E(\varphi_{s_2},\psi_{s_2};s_2)+|s_1-s_2|E_1(\varphi_{s_2},\psi_{s_2})\\[2mm]
&=\mu(s_2)+|s_1-s_2|E_1(\varphi_{s_2},\psi_{s_2}).
\end{array}
\end{eqnarray*}
Putting these two inequalities together and employing (\ref{2.2.0117}), we conclude that
\begin{eqnarray*}  %\label{2.2.00110}
\mu(s_1)\leq \mu(s_2)+K|s_1-s_2|.
\end{eqnarray*}

Reversing the role of the indices 1 and 2 in the derivation of the above inequality gives
the same bound with the indices switched. Therefore, we deduce that
\begin{eqnarray*}    %%%% \label{2.2.00111}
|\mu(s_1)-\mu(s_2)|\leq K|s_1-s_2|,
\end{eqnarray*}
which completes the proof of Proposition 2.4. \hfill $\Box$
\vspace{2mm}

To emphasize the dependence on the parameters, we write
\begin{eqnarray*} %\label{2.2.111}
\varphi=\varphi_s(|\xi|,x_3),\ \psi=\psi_s(|\xi|,x_3), \mbox{ and }\lambda=\lambda(|\xi|,s).
\end{eqnarray*}
In view of Propositions \ref{pro:0203} and \ref{pro:0201}, we can state the following
existence result of solutions to (\ref{2.2.4}),
(\ref{2.2.5}) and (\ref{2.2.3}) for each $(\xi, s)\in \mathcal{M}$.

\begin{Proposition} For each $(\xi, s)\in \mathcal{M}$, there exists a solution
$\varphi_s(|\xi|,x_3)$, $\psi_s(|\xi|,x_3)$ with
$\lambda=\lambda(|\xi|,s)>0$ to the problem (\ref{2.2.4}) along with the corresponding
jump and boundary conditions. Moreover, these solutions $\psi_s(|\xi|,0)\neq 0$ and
the solutions are smooth when restricted to either $(-m,0)$ or $(0,l)$.
\end{Proposition}

In order to prove $\lambda\rightarrow +\infty$ as $|\xi|\rightarrow
+\infty$, and give the existence of solutions to the original
problem (\ref{2.2.040}), (\ref{2.2.04}) and (\ref{2.2.3}), we shall
further restrict  $\xi$ to satisfy
\begin{eqnarray*}
|\xi|> R_2:=2C^{-1}_0\sqrt{C_2}\geq R_1, \mbox{ and }C_2=\max\{4C_1,(R_1
C_0/2)^2\},\end{eqnarray*} where $C_0$, $C_1$  and  $R_1$ are the
constants from Proposition 2.1. Thus we have the following
conclusion.

\begin{Theorem}
For each $\xi$ with $|\xi|> R_2\geq R_1$, there exists a $s\in (0,s_1):=(0,C_0R_2/C_2)$,
such that
\begin{equation}\label{2.2.12}   s=1/\lambda^2(|\xi|,s).   \end{equation}
\end{Theorem}
\begin{Remark} It is easy to check that $(\xi,s)\in \mathcal{M}$ for
any $|\xi|>R_2$ and $s\in (0,s_1)$.
\end{Remark}
{\bf{Proof.}} Recalling $-\lambda^2(s)=\mu(s)$, we define
\begin{eqnarray*} % \label{2.2.13}
F(s)=s\lambda^2(s)-1=-s\mu(s)-1.
\end{eqnarray*}
According to Proposition 2.4, $F(s)$ is continuous. Moreover, we have
\begin{eqnarray}\label{2.2.14}
F(0)=-1<0.
\end{eqnarray}
Now if  there exists a $\bar{s}\in(0,s_1)$ such that
\begin{eqnarray}\label{2.2.15}
F(\bar{s})>0,
\end{eqnarray}
then combining (\ref{2.2.14}) with (\ref{2.2.15}), we can find $s\in (0,\bar{s})$
 such that $F(s)=0$, i.e., (\ref{2.2.12}) holds. In the following, we verify (\ref{2.2.15}).

According to (\ref{2.2.10}) and the fact $C_2> C_1$, we see that
\begin{eqnarray*} % \label{2.2.16}
-\mu(s)> C_0R_2-sC_1\geq  C_0R_2-sC_2,
\end{eqnarray*}
which yields
\begin{eqnarray*} % \label{2.2.17}
F(s)> -C_2s^2+C_0R_2s-1:=f(s).
\end{eqnarray*}
Recalling the definition of $s_1$ and $R_2$, it is easy to verify that
there exists at least one positive root $0<\bar{s}<s_1$, such that
$f(\bar{s})=0$. This implies (\ref{2.2.15}).  \hfill $\Box$

\subsection{Construction of a solution to (\ref{2.2.1})}

We may now use Theorem 2.1 to think of $s=s(|\xi|)$, since we can find $s\in {\mathcal{S}}$
so that (\ref{2.2.12}) holds. As such we may also write $\lambda=\lambda(|\xi|)$ from now on.

Once Theorem 2.1 is established, we can combine it with Proposition 2.5 to obtain immediately
a solution to (\ref{2.2.04}), and in turn a solution to (\ref{2.2.1}) for each spacial
frequency $\xi$ with $|\xi|>R_2$.
%%%%%%%%%%%%%%%%%%%%%%%%%%%%%%%%%%%%%%%%%%%%%%
\begin{Theorem} For $\xi\in {\mathbb{R}}^2$ with $|\xi|>R_2$, there exists a solution
$\varphi=\varphi(\xi,x_3)$, $\theta=\theta(\xi,x_3)$, $\psi=\psi(\xi,x_3)$
 and $\lambda=\lambda(|\xi|)>0$ to (\ref{2.2.1}), such that $\psi(0)\neq 0$.
 The solution is smooth when restricted to $(-m,0)$ or $(0,l)$,
 and it is equivariant in $\xi$ in the sense that if $\mathcal{R}\in \mathcal{SO}(2)$
 is a rotation operator, then
 \begin{eqnarray*} % \label{2.2.18}
 \left(\begin{array}{ll}
 \varphi(\mathcal{R}\xi,x_3)\\[2mm]
 \theta(\mathcal{R}\xi,x_3)\\[2mm]
 \psi(\mathcal{R}\xi,x_3)
 \end{array}
 \right)=
 \left(
 \begin{array}{ll}
 \mathcal{R}_{11}\ \mathcal{R}_{12}\ 0\\[2mm]
 \mathcal{R}_{21}\ \mathcal{R}_{22}\ 0\\[2mm]
 0\ \ \ \ 0\ \ \ \ 1
 \end{array}
 \right)
 \left(
 \begin{array}{ll}
 \varphi(\xi,x_3)\\[2mm]
 \theta(\xi,x_3)\\[2mm]
 \psi(\xi,x_3)
 \end{array}
 \right)
 \end{eqnarray*}
Moreover,
\begin{eqnarray}\label{2.2.010}
\lambda^2\geq C_3|\xi|-1:=C_1|\xi|/2-1\ \mbox{ for }\  |\xi|\geq
\max\left\{\frac{C_1+1}{C_0},R_2\right\},
\end{eqnarray} where $C_0$, $C_1$  and $R_2$ are the constants from
Proposition 2.1 and Theorem 2.1.
\end{Theorem}
{\bf{Proof}}. In view of  Theorem 2.1 and Proposition 2.5, there
exists a solution ($\varphi(|\xi|),\psi(|\xi|),\theta\equiv0)$ to
(\ref{2.2.1}) for the fixed frequency given by $(|\xi|,0)$ with
$|\xi|>R_2$. Thus we may find a rotation operator $\mathcal{R}\in
\mathcal{SO}(2) $ so that $\mathcal{R}(\xi)=(|\xi|,0)$. Define $(\varphi(\xi,x_3),
\theta(\xi,x_3))=\mathcal{R}^{-1}(\varphi(|\xi|,x_3),0)$ and
$\psi(\xi,x_3)=\psi(|\xi|,x_3)$. This gives a solution to
$(\ref{2.2.1})$  for any frequency $\xi$ with $|\xi|>R_2$.
The equivalence in $\xi$ follows from the definition.

We proceed to estimate (\ref{2.2.10}). According to (\ref{2.2.9}), for $|\xi|> R_2\geq R_1$,
one has
\begin{eqnarray}\label{2.2.00111}
-\lambda^2\leq -C_0|\xi|+sC_1.
\end{eqnarray}
The fact that $s=1/\lambda^2$ from Theorem 2.1 combined with
(\ref{2.2.00111}) results in
\begin{eqnarray*}  %\label{2.2.001112}
\lambda^4-C_0|\xi|\lambda^2+C_1\geq 0,
\end{eqnarray*}
which leads to
\begin{eqnarray}\label{2.2.01113}
\lambda^2\geq
\frac12\left(C_0|\xi|+\sqrt{C_0^2|\xi|^2-4C_1}\right)>\frac12\sqrt{C_0^2|\xi|^2-4C_1},
\end{eqnarray}or
\begin{eqnarray}\label{5.5.01113}
\lambda^2\leq
\frac12\left(C_0|\xi|-\sqrt{C_0^2|\xi|^2-4C_1}\right)<\sqrt{C_1}.
\end{eqnarray}

On the other hand, in view of Theorem 2.1, we find that
\begin{eqnarray*}
\lambda^2>C_0R_2/C_2=C_2^{1\over 2}/2\geq \sqrt{C_1},
\end{eqnarray*}
which implies that $(\ref{5.5.01113})$ does not hold. Hence,
$\lambda^2$ has to satisfy (\ref{2.2.01113}).

 For $|\xi|>{(C_1+1)}/{C_0}$, we have
\begin{eqnarray}\label{2.2.01114}
\sqrt{C_0^2|\xi|^2-4C_1}>\sqrt{C_1^2|\xi|^2-4C_1|\xi|+4}=\sqrt{(C_1|\xi|-2)^2}=C_1|\xi|-2.
\end{eqnarray}
By virtue of (\ref{2.2.01113}) and (\ref{2.2.01114}), the estimate (\ref{2.2.010})
is obtained by denoting $C_3=C_1/2$. \hfill $\Box$
\vspace{2mm}

Next, we derive an estimate for the $H^k$-norm of the
solutions $(\varphi, \psi)$ with $|\xi|$ varying, which will be used in the proof of
Theorem 2.3 when integrating solutions in a Fourier synthesis.
\begin{Lemma}
Let $\varphi(|\xi|)$, $\psi(|\xi|)$ be the solutions to
(\ref{2.2.4}) constructed in Theorem 2.2. Let $R_2$, $C_0$, $C_1$,
$C_3>0$ be the constants from Theorem 2.2. Fix
$R_0=\max\{{(C_1+1)}/{C_0}, 2/C_3,R_2,1\}\in(0,\infty)$. Then for
any $\xi$ with $|\xi|\geq R_0$ and for each $k\geq 0$, there exists a
constant $A_k>0$ depending on $\rho_0$, $p$, $g$, $\omega$, $l$ and $m$, such that
\begin{eqnarray}\label{2.2.19}
\|\varphi(|\xi|)\|_{H^k{(-m,0)}}+\|\psi(|\xi|)\|_{H^k{(-m,0)}}+\|\varphi(|\xi|)
\|_{H^k{(0,l)}}+\|\psi(|\xi|)\|_{H^k{(0,l)}}\leq A_k\sum_{j=0}^k|\xi|^j.
\end{eqnarray}
Also, there exists a constant $B_0>0$ depending on the same parameters,
such that for any $|\xi|>0$,
\begin{eqnarray}\label{2.2.19a}
\left\|\sqrt{\varphi^2(|\xi|)+\psi^2(|\xi|)}\right\|_{L^2(-m,l)}\geq B_0.
\end{eqnarray}
\end{Lemma}
{\bf{Proof.}} We begin with the proof of (\ref{2.2.19}). For
simplicity we will derive an estimate of the $H^k$-norm on the
interval $(0,l)$ only. A bound on $(-m,0)$ follows similarly, and the desired
result follows from adding the two estimates together. First note that the choice of $R_0$,
when combined with Theorem 2.2 and (\ref{2.3.512223}), implies that
\begin{eqnarray}\label{2.2.20}
-g|\xi|\leq \mu \leq 1-C_3|\xi|=1-\frac{C_3}{2}|\xi|-\frac{C_3}{2}|\xi|\leq -\frac{C_3}{2}|\xi|.
\end{eqnarray}
Also keep in mind that $\rho_0$ is smooth on each interval $(0,l)$ and
$(-m,0)$ and bounded from above and below. Throughout the proof we denote by $C$
a generic positive constant depending on the appropriate parameters.
We proceed by induction on $k$. For $k=0$ the fact that
$(\varphi(|\xi|), \psi(|\xi|))\in {\mathcal{A}}$ implies that there
is a constant $A_0>0$ depending on the various parameters, so that
\begin{eqnarray*}  %%\label{2.2.21}
\|\varphi(|\xi|)\|_{L^2{(0,l)}}+\|\psi(|\xi|)\|_{L^2{(0,l)}}\leq A_0.
\end{eqnarray*}
Suppose now that the bound holds for some $k-1\geq 0$, i.e.,
\begin{eqnarray*}  %%\label{2.2.22}
\|\varphi(|\xi|)\|_{H^{k-1}{(0,l)}}+\|\psi(|\xi|)\|_{H^{k-1}{(0,l)}}\leq
A_{k-1}\sum_{j=0}^{k-1}|\xi|^j.
\end{eqnarray*}
Define $w(|\xi|):=
p'_+(\rho_0)\rho_0(\psi'(|\xi|)+|\xi|\varphi(|\xi|))$, where $'=\partial_{x_3}$.
Recalling $\mu=-\lambda^2$ and $s=\lambda^{-2}$, the system (\ref{2.2.4}) implies that
\begin{eqnarray}\label{2.2.23}
\left\{
\begin{array}{ll}
w(|\xi|) =\mu|\xi|^{-1}\rho_0\varphi(|\xi|)+g\rho_0\psi(|\xi|)
+{4\omega^2\rho_0}({\mu}|\xi|)^{-1}\varphi(|\xi|) \\[2mm]
w'(|\xi|)=-\mu\rho_0\psi-g|\xi|\rho_0\varphi.
\end{array}
\right.
\end{eqnarray}
These equations, together with (\ref{2.2.20}) and the fact that $|\xi|>1$, yield that
\begin{eqnarray*}  %%\label{2.2.24}
\begin{array}{ll}
\|w(|\xi|)\|_{H^{k-1}(0,l)}\leq C\left(\|\varphi(|\xi|)\|_{H^{k-1}(0,l)}+\|\psi(|\xi|)\|_{H^{k-1}(0,l)}\right),
\\[3mm]
\|w'(|\xi|)\|_{H^{k-1}(0,l)}\leq
C|\xi|\left(\|\varphi(|\xi|)\|_{H^{k-1}(0,l)}+\|\psi(|\xi|)\|_{H^{k-1}(0,l)}\right),
\end{array}
\end{eqnarray*}
so that $\|w(|\xi|)\|_{H^{k}(0,l)}\leq C \sum_{j=0}^{k}|\xi|^j$.

On the other hand, the definition of $w(|\xi|)$ implies that
\begin{eqnarray*}  %%\label{2.2.25}
\|\psi'(|\xi|)\|_{H^{k-1}(0,l)}\leq
\left\|\frac{w(|\xi|)}{p'(\rho_0)\rho_0}\right\|_{H^{k-1}(0,l)}+|\xi|\|\varphi(|\xi|)\|_{H^{k-1}(0,l)},
\end{eqnarray*}
such that $\|\psi(|\xi|)\|_{H^{k}(0,1)}\leq C \sum_{j=0}^{k}|\xi|^j$.
Returning to the first equation in (\ref{2.2.23}), we see that
\begin{eqnarray*}
\varphi(|\xi|)\leq\frac{|\xi|[w(|\xi|)-g\rho_0\psi(|\xi|)]}{\mu+4\rho_0\omega^2/\mu}
=\frac{|\xi|}{\mu}\frac{[w(|\xi|)-g\rho_0\psi(|\xi|)]}{1+4\rho_0\omega^2/\mu^2},
\end{eqnarray*}
from which we get
\begin{eqnarray*}  %%\label{2.2.26}
\|\varphi(|\xi|)\|_{H^{k}(0,l)}\leq
\frac{|\xi|}{|\mu|}\left[g\|\rho_0\psi(|\xi|)\|_{H^{k}(0,l)}+\|w(|\xi|)\|_{H^{k}(0,l)}\right].
\end{eqnarray*}
Thus, making use of (\ref{2.2.20}) we deduce that
$\|\varphi(|\xi|)\|_{H^{k}(0,l)}\leq C \sum_{j=0}^{k}|\xi|^j$. Hence,
\begin{eqnarray*}  %%\label{2.2.27}
\|\varphi(|\xi|)\|_{H^{k}(0,l)}+\|\varphi(|\xi|)\|_{H^{k}(0,l)}\leq
A_{k}\sum_{j=0}^{k}|\xi|^j
\end{eqnarray*}
for some constant $A_{k}>0$ depending on the parameters, i.e., the bound holds for $k$.
By induction, the bound holds for all $k\geq 0$. Finally,
(\ref{2.2.19a}) follows from the fact that $(\varphi(|\xi|), \psi(|\xi|))\in {\mathcal{A}}$
and $\rho_0$ is bounded from above and below. \hfill $\Box$

\subsection{Fourier synthesis}

 Inspired by \cite{Y-I}, we will use the Fourier synthesis to build growing solutions to
(\ref{1.2.14}) out of the solutions constructed in the previous
section for fixed spacial frequency $\xi\in\mathbb{R}^2$. The
solutions will be constructed to grow in the piecewise sobolev space
of order $k$, $H^k$, defined by (\ref{1.4.3}).
%%%%%%%%%%%%%%%%%%%%%%%%%%%%%%%%%%%%%%%%%%%%%%%%%%%%%%%%%%%%%%
\begin{Theorem}
Let $R_0\leq R_3<R_4<\infty$, where $R_0>1$. Let $f\in C_c^\infty(\mathbb{R}^2)$
be a real-valued function, such that
$f(\xi)=f(|\xi|)$ and ${\rm{supp}}(f)\subset {B(0, R_4)/B(0,R_3)}$.
For $\xi\in \mathbb{R}^2$ define
\begin{eqnarray*}
\begin{array}{ll}
\hat{w}(\xi,x_3)=-i\varphi(\xi,x_3)e_1-i\theta(\xi,x_3)e_2+\psi(\xi,x_3)e_3,
\end{array}
\end{eqnarray*}
where $\varphi,\theta,\psi$ and $\pi$ are solutions to (\ref{2.2.1}).
Writing $x'\cdot\xi=x_1\xi_1+x_2\xi_2$, we define
\begin{eqnarray}\label{2.4.1}
\left\{\begin{aligned}
&\eta(t,x)=\frac{1}{4\pi^2}\displaystyle\int_{\mathbb{R}^2}f(\xi)\hat{w}(\xi,x_3)e^{\lambda(|\xi|)t}e^{ix'\cdot
\xi}d\xi,\\
&v(t,x)=\frac{1}{4\pi^2}\displaystyle\int_{\mathbb{R}^2}\lambda(|\xi|)f(\xi)\hat{w}(\xi,x_3)
e^{\lambda(|\xi|)t}e^{ix'\cdot  \xi}d\xi,\\
&q(t,x)=-\frac{\rho_0(x_3)}{4\pi^2}\displaystyle\int_{\mathbb{R}^2}f(\xi)
(\xi_1\varphi(\xi,x_3)+\xi_2\theta(\xi,x_3)+\partial_{x_3}\psi(\xi,x_3))
e^{\lambda(|\xi|)t}e^{ix'\cdot \xi}d\xi,
\end{aligned}\right.
\end{eqnarray} where we have defined
$\lambda(|\xi|)=\varphi(\xi,x_3)\equiv\theta(\xi,x_3)\equiv\psi(\xi,x_3)\equiv0$
if $|\xi|\leq R_0$. Then $\eta,\  v,\ q$ are real-valued solutions to the
linearized equations (\ref{1.2.14}) along with the corresponding jump
and boundary conditions. For every $k\in\mathbb{N}$, we have the estimate
\begin{eqnarray}\label{2.4.2}
\|\eta(0)\|_{H^k}+\|v(0)\|_{H^k}+\|q(0)\|_{H^k}\leq
\bar{C}_k\left(\displaystyle\int_{\mathbb{R}^2}(1+|\xi|^2)^{k+1}|f(\xi)|^2d\xi\right)^{1/2}<\infty
\end{eqnarray}
for some constant $\bar{C}_k>0$ depending on the parameters $\rho_0$,
$p$, $l$ $m$, $\omega$ and $g$. Furthermore, for every $t>0$ we have
$\eta(t),v(t),q(t)\in H^k(\Omega)$ and
\begin{eqnarray}\label{2.4.3}
\left\{
\begin{array}{ll}
e^{t\sqrt{C_4R_3-1}}\|\eta(0)\|_{H^k}\leq \|\eta(t)\|_{H^k}\leq
e^{t\sqrt{ gR_4}}\|\eta(0)\|_{H^k},\\[2mm]
e^{t\sqrt{C_4R_3-1}}\|v(0)\|_{H^k}\leq \|v(t)\|_{H^k}\leq
e^{t\sqrt{ gR_4}}\|v(0)\|_{H^k},\\[2mm]
e^{t\sqrt{C_4R_3-1}}\|q(0)\|_{H^k}\leq \|q(t)\|_{H^k}\leq
e^{t\sqrt{ gR_4}}\|q(0)\|_{H^k} .
\end{array}
\right.
\end{eqnarray}
\end{Theorem}
{\bf{Proof}}. For each fixed $\xi\in\mathbb{R}^2$,
\begin{eqnarray*} &&
\eta(x,t)=f(\xi)\hat{w}(\xi,x_3)e^{\lambda(|\xi|)t}e^{ix'\cdot\xi} , \\[1mm]
&& v(x,t)=\lambda(|\xi|)f(\xi)\hat{w}(\xi,x_3)e^{\lambda(|\xi|)t}e^{ix'\cdot\xi}, \\[1mm]
&& q(x,t)=-\rho_0(x_3)f(\xi)(\xi_1\varphi(\xi,x_3)+\xi_2\theta(\xi,x_3)
+\partial_3\psi(\xi,x_3))e^{\lambda(|\xi|)t}e^{ix'\cdot\xi}
\end{eqnarray*}
gives a solution to (\ref{1.2.14}). Since ${\rm{supp}}(f)\subset B(0,R_4)/B(0,R_3)$,
Lemma 2.2 results in
\begin{eqnarray*}
\sup\limits_{\xi\in {\rm{supp}}(f)}\|\partial_{x^3}^k
\hat{w}(\xi,\cdot)\|_{L^\infty} <\infty \mbox{ for all }k\in \mathbb{N}.
\end{eqnarray*}
These bounds imply that the Fourier synthesis of the solutions given by (\ref{2.4.1})
is also a solution to (\ref{1.2.14}). The bound (\ref{2.4.2}) follows from Lemma 2.2
with arbitrary $k\geq 0$ and the fact that $f$ is compactly supported. From
(\ref{2.2.20}), the estimates (\ref{2.4.3}) follow. \hfill $\Box$

%%%%%%%%%%%%%%%%%%%%%%%%%%%%%%%%%%%%%%%%%%%%%%%%%%%%%%%%%
\section{Ill-posedness for the linear problem}
\quad \ We assume that $\eta$, $v$, $q$ are the real-valued solutions to (\ref{1.2.14})
along with the corresponding jump and boundary conditions established in Theorem 1.1.
Furthermore, suppose that the solutions are band-limited at radius
$R>0$, i.e. that
\begin{eqnarray*}  %%%\label{0401}
\bigcup_{x_3\in (-m,l)}{\mathrm{supp}}(|\hat{\eta}(\cdot,x_3)|
+|\hat{v}(\cdot, x_3)|+|\hat{q}(\cdot, x_3)|)\subset B(0,R),\end{eqnarray*}
where $\hat{v}$ denotes the horizontal Fourier transform defined by (\ref{0113}).
We will derive estimates for band-limited solutions in terms of $R$.

Differentiating the second equation in (\ref{1.2.14}) with respect to times $t$
 and eliminating the $\eta$ term by using the first equation, we obtain
\begin{equation} \label{0402}
\rho_0\partial_{tt}v-\nabla(p'(\rho_0)\rho_0{\mathrm{div}}
v)+g\rho_0{\mathrm{div}} v e_3+2\rho_0\omega \partial_tv_{2}
e_1-2\rho_0\omega
\partial_t v_{1} e_2=0
\end{equation}
along with the jump and boundary conditions
\begin{eqnarray*} && \llbracket \partial_tv_3 \rrbracket =0,
\ {\mbox{ and }} \ \llbracket p'(\rho_0)\rho_0 {\mathrm{div}} v\rrbracket=0, \\[1mm]
&& \partial_t v_3(t,x',-m)=\partial_tv_{3}(t,x',l)=0.
\end{eqnarray*}
The band limited assumption implies that
supp$(\hat{v}(\cdot,x_3))\subset B(0,R)$ for all $x_3\in (-m,l)$.
The initial datum for $\partial_t v(0)$ is given in terms of the
initial data $q(0)$ and $\eta(0)$ via the second linear equations,
i.e.,
\begin{eqnarray*}   %%%\label{0405}
\rho \partial_t v(0)=-\nabla
(p'(\rho_0)q(0))-gq(0)e_3-g\rho_0\nabla\eta_3(0)+2\rho_0\omega
v_2(0)e_1-2\rho_0\omega v_1(0)e_2.
\end{eqnarray*}

By the standard energy estimate procedure (see \cite[Section 3]{Y-I2}),
we can establish three energy estimates in the following which
ensure the uniqueness result in Theorem 3.1.

\begin{Lemma} For solutions to (\ref{0402}) it holds that
\begin{eqnarray*}  %%%\label{0406}
\partial_t\displaystyle\int_{\Omega}\left(\frac{\rho_0}{2}|\partial_t v|^2
+\frac{p'(\rho_0)\rho_0}{2}\left|{\mathrm{div}}v-\frac{g}{p'(\rho_0)}v_3\right|^2\right)
=\partial_t\displaystyle\int_{\mathbb{R}^2}
\frac{g\llbracket \rho_0 \rrbracket}{2}|v_3|^2.
\end{eqnarray*}
\end{Lemma}

\begin{Lemma}
Let $v\in H^1(\Omega)$ be band-limited at radius $R>0$ and satisfy
the boundary conditions $v_3(t,x',-m)=v_3(t,x',l)=0$. Then
\begin{eqnarray*}  %%%\label{0407}
\displaystyle\int_{\Omega}\frac{g\llbracket \rho_0
\rrbracket}{2}|v_3|^2-\frac{p'(\rho_0)\rho_0}{2}\left|{\mathrm{div}}v-\frac{g}{p'(\rho_0)}v_3\right|^2\leq
\frac{\Lambda^2({R})}{2}\int_{\Omega}\rho_0|v|^2,
\end{eqnarray*}
where the value of $\Lambda(R)$, depending $R$, is given by (4.2) in \cite[Section 4.1]{Y-I}.
\end{Lemma}

\begin{Proposition}
Let $v$ be a solution to (\ref{0402}) along with the corresponding jump and boundary
conditions that is also band-limited at radius $R>0$. Then
$$
\|v(t)\|^2_{L^2(\Omega)}+\|\partial_t v(t)\|^2_{L^2(\Omega)} \leq Ce^{2\Lambda(R)t}
\left(\|v(0)\|^2_{L^2(\Omega)} +\|\partial_t v(0)\|^2_{L^2(\Omega)}
+\|{\mathrm{div}}v(0)\|^2_{L^2(\Omega)}\right) $$
for a constant $C=C(\rho_0,l,m,p,g,\Lambda(R))>0$, where the value of
$\Lambda(R)$, depending $R$, is given by (4.2) in \cite[Section 4.1]{Y-I}.
\end{Proposition}

 Similar to \cite{Y-I}, once we get Proposition 3.1, through
constructing the horizontal spatial frequency projection operator,
we can obtain the uniqueness result. Here we give the proof for the
reader's convenience.

Let $\Phi\in C^\infty_0(\mathbb{R}^2)$ be so that $0\leq \Phi \leq
1$, ${\rm{supp}}(\Phi)\subset B(0,1)$, and $\Phi(x)=1$ for $x\in
B(0,1/2)$. For $R>0$ let $\Phi_R$ be the function defined by
$\Phi_R(x)=\Phi(x/R)$. We define the projection operator $P_R$ via
\begin{eqnarray*}  %%%\label{0422}
P_Rf={\mathcal{F}}^{-1}(\Phi_R{\mathcal{F}}f),\qquad f\in L^2(\Omega),
\end{eqnarray*}
where $\mathcal{F}\cdot =\hat{\cdot}$ denotes the horizontal Fourier transform
in $x'$. It is easy to see that $P_R$ satisfies the following.
\begin{enumerate}
  \item[(1)] $P_R f$ is  band-limited at radius $R$.
  \item[(2)]  $P_R$ is a bounded linear operator on $H^k(\Omega)$ for all $k\geq 0$.
  \item[(3)] $P_R$ commutes with partial differentiation and
multiplication by functions depending only on $x_3$.
  \item[(4)] $P_R f=0$ for all $R>0$ if and only if  $f=0$.
\end{enumerate}

Now we are able to prove the uniqueness results on $\eta$, $v$ and $q$.
\begin{Theorem}
Assume that $(\eta_1,v_1,q_1)$ and $(\eta_2,v_2,q_2)$  are two
solutions to (\ref{1.2.14}). Then $\eta_1=\eta_2$, $v_1=v_2$ and
$q_1=q_2$.
\end{Theorem}
{\bf{Proof}}.
It suffices to show that solutions to (\ref{1.2.14}) with $0$ initial
data remain $0$ for $t>0$. Suppose that $\eta$, $v$ are
solutions with vanishing initial data. Fix $R>0$ and define
$\eta_R=P_R\eta$, $v_R=P_Rv$, $q_R=P_Rq$. The properties of $P_R$
show that $\eta_R$, $v_R$, $q_R$ are also solutions to (\ref{1.2.14})
but that they are band-limited at radius $R$. Turning to the second
order formulation, we find that $v_R$ is a solution to (\ref{0402})
with initial data $v_R(0)=\partial_tv_R(0)=0$. We may then apply
Proposition 3.1 to deduce that
\begin{eqnarray*}  %%%\label{0423}
\|v_R(t)\|_{L^2(\Omega)}=\|\partial_t v_R(t)\|_{L^2(\Omega)}=0\mbox{
for all }t\geq 0,
\end{eqnarray*}
which shows that $\eta_R(t)$, $v_R(t)$ and $q_R(t)$ all vanish for $t\geq 0$.
Since $R$ is arbitrary, it must hold that $\eta(t)$, $v(t)$ and $q(t)$
 also vanish for $t\geq 0$. \hfill $\Box$
 \vspace{2mm}

The solutions to the linear problem (\ref{1.2.14}) constructed in
Theorem 2.3 are sufficiently pathological to give rise to a result
showing that the solutions depend discontinuously on the initial
data. Then, in spite of the previous uniqueness result, we obtain that the
linear problem is ill-posed in the sense of Hadamard. Next, we prove Theorem 1.1.

 Similar to \cite{Y-I}, once we get the pathological solutions to the linear problem
 (\ref{1.2.14}) constructed in Theorem 2.3, we are able to  show that
the solutions depend discontinuously on the initial data described in Theorem 1.1.
 Here we give the proof for the reader's convenience.
\\[0.5em]
{\bf{Proof of Theorem 1.1}}. Fix $j\geq k\geq 0$, $\alpha>0$,
$T_0>0$ and let $\bar{C}_j$, $R_0$, $B_0$, $C_4>0$ be the constants
from Theorem 2.3, Lemma 2.2 and Theorem 2.2 respectively. For each
$n\in \mathbb{N}$, let $R(n)$ be sufficiently large so that
$R(n)>R_0$, $\sqrt{C_4R(n)-1}\geq 1$, and
\begin{eqnarray*}  %%\label{3.3.3}
\frac{\exp(2T_0\sqrt{C_4R(n)-1})}{(1+(R(n)+1)^2)^{j-k+1}(1+R(n))^2}\geq \alpha^2 n^2\bar{C}_j^2/B_0^2.
\end{eqnarray*}
Choose $f_n\in C_c^\infty(R^2)$, so that ${\rm{supp}}(f_n)\subset
B(0,R(n)+1)/B(0,R(n))$, $f_n$ is real-valued and radial, and
\begin{eqnarray}\label{3.3.4}
\displaystyle\int_{\mathbb{R}^2}(1+|\xi|^2)^{j+1}
|f_n(|\xi|)|^2d\xi=\frac{1}{\bar{C}_j^2 n^2}.
\end{eqnarray}

We may now apply Theorem 2.3 with $f_n,\ R_3=R(n)$, and $R_4=R(n)+1$ to find
$\eta_n, \ v_n,\ q_n$ that solve (\ref{1.2.14}) with the corresponding jump
and boundary conditions, such that $\eta_n, \ v_n,\ q_n\in H^j(\Omega)$ for all $t\geq 0$.
By (\ref{2.4.2}) and the choice of $f_n$ satisfying (\ref{3.3.4}), we find that (\ref{3.3.1})
holds for all $n$.

On the other hand, we have
\begin{equation*}  %%%\label{3.3.5}
\begin{aligned}
\|\eta_n(T_0)\|_{H^k}^2&\geq \displaystyle\int_{\mathbb{R}^2}(1+|\xi|^2)^k
|f_n(\xi)|^2e^{2T_0\lambda(|\xi|)}\|\hat{w}(\xi,\cdot)\|_{L^2(-m,l)}^2 d\xi\\
&\geq\frac{\exp(2T_0\sqrt{C_4R(n)-1})}{(1+(R(n)+1)^2)^{j-k+1}}
\displaystyle\int_{\mathbb{R}^2}(1+|\xi|^2)^{j+1}|f_n(\xi)|^2\|\hat{w}(\xi,\cdot)\|_{L^2 (-m,l)}^2 d\xi\\
&\geq\frac{\alpha^2 n^2\bar{C}_j^2}{B_0^2}\displaystyle\int_{\mathbb{R}^2}(1+|\xi|^2)^k|f_n(\xi)|^2
B_0^2 d\xi=\alpha^2.
\end{aligned}
\end{equation*}
Here the first bound is trivial, while the second one follows from ${\rm{supp}}(f_n)\subset B(0,R(n)+1)$
and $\lambda(|\xi|)\geq \sqrt{C_4R(n)-1}$, and the third one from the choice of $R(n)$ and
the lower bound (\ref{2.2.19a}). Since $\lambda(|\xi|)\geq \sqrt{C_4R(n)-1}\geq 1$
on the support of $f_n$, we deduce that
\begin{eqnarray*}  %%\label{3.3.6}
\|v_n(t)\|_{H^k}^2\geq \|\eta_n(t)\|_{H^k}^2\geq\|\eta_n(T_0)\|_{H^k}^2 \  {\rm{for}} \ {\rm{all}}\ t\geq 0.
\end{eqnarray*}
from which (\ref{3.3.2}) follows. \hfill $\Box$

\section{Proof of Theorem 1.2}
\setcounter{equation}{0}
\quad \ The proof is similar to \cite{Y-I}
under necessary modifications. We argue by contradiction. Suppose
 that the perturbed problem has property $EE(k)$ for some $k\geq 3$.
 Let $\delta,\ t_0$ and $C>0$ be the constants and function provided by the property $EE(k)$.
 Fix $n\in \mathbb{N}$ so that $n>C$, applying Theorem 1.1
with this $n$, $T_0=t_0/2$, $k\geq 3$ and $\alpha=1$, we can find
$\bar{\eta}$, $\bar{v}$, $\bar{\sigma}$ solving (\ref{1.2.14}), such that
\begin{eqnarray*}  %%\label{4.1.12}
\|(\bar{\eta}, \bar{v},\bar{\sigma})(0)\|_{H^k}<\frac 1n,
\end{eqnarray*}
but
\begin{eqnarray}  \label{4.1.12a}
\|\bar{v}(t)\|_{H^3}\geq\|\bar{\eta}(t)\|_{H^3}\geq 1\ {\rm{for}}\ t\geq t_0/2.
\end{eqnarray}
For $\epsilon>0$ we define $\bar{\eta}_0^\epsilon=\epsilon\bar{\eta}(0),\
\bar{v}_0^\epsilon=\epsilon\bar{v}(0),$ and $\bar{\sigma}_0^\epsilon=\epsilon\bar{\sigma}(0)$.

Then for $\epsilon<\delta n$, we have $\|(\bar{\eta}^\epsilon_0,
\bar{v}^\epsilon_0, \bar{\sigma}^\epsilon_0)\|_{H^k}<\delta$. So,
according to $EE(k)$ there exist $\tilde{\eta}^\epsilon,
v^\epsilon,\sigma^\epsilon\in L^\infty(0,t_0;H^3(\Omega))$ that
solve the perturbed problem with $(\bar{\eta}_0^\epsilon,
\bar{v}_0^\epsilon, \bar{\sigma}_0^\epsilon)$ as initial data and
that satisfy the inequality
\begin{eqnarray} \label{4.1.13}
\sup\limits_{0\leq t<t_0}\|(\tilde{\eta}^\epsilon, v^\epsilon,
\sigma^\epsilon)(t)\|_{H^3}\leq F(\|(\bar{\eta}_0^\epsilon,
\bar{v}_0^\epsilon, \bar{\sigma}_0^\epsilon)\|_{H^k})\leq
C\epsilon\|(\bar{\eta}, \bar{v}, \bar{\sigma})(0)\|_{H^k}<\epsilon.
\end{eqnarray}

Now, defining the rescaled functions
$\bar{\eta}^\epsilon=\tilde{\eta}^\epsilon/\epsilon$,
$\bar{v}^\epsilon=v^\epsilon/\epsilon$, $\bar{\sigma}^\epsilon
=\sigma^\epsilon/\epsilon$, and rescaling (\ref{4.1.13}), we infer that
\begin{eqnarray}\label{4.1.14}
\sup\limits_{0\leq t<t_0}\|(\bar{\eta}^\epsilon, \bar{v}^\epsilon,
\bar{\sigma}^\epsilon)(t)\|_{H^3}\leq 1.
\end{eqnarray}
On the other hand, $(\bar{\eta}^\epsilon, \bar{v}^\epsilon,\bar{\sigma}^\epsilon)$ satisfies
\begin{eqnarray}\label{4.141}
\left\{
\begin{array}{ll}
\partial_t\bar{\eta}^\epsilon=\bar{v}^\epsilon,\\[3mm]
\partial_t\bar{\sigma}^\epsilon+\rho_0 {\mathrm{div}}{\bar{v}^\epsilon}+\epsilon(\bar{\sigma}^\epsilon
{\mathrm{div}}\bar{v}^\epsilon-\rho_0 {\mathrm{tr}}
(\bar{B}^\epsilon D\bar{v}^\epsilon))-
\epsilon^2(\bar{\sigma}^\epsilon {\mathrm{tr}}(\bar{B}^\epsilon D\bar{v}^\epsilon))=0,\\[2mm]
\left(\partial_t\bar{v}^\epsilon+
\nabla(h'(\rho_0)\bar{\sigma}^\epsilon+g
e_3\cdot\bar{\eta}^\epsilon\right)+
\nabla\bar{\mathfrak{R}}^\epsilon
-\epsilon(\bar{B}^\epsilon\nabla(h'(\rho_0)\bar{\sigma}^\epsilon+g\bar{\eta}^\epsilon)
+\bar{B}^\epsilon\nabla\bar{\mathfrak{R}}^\epsilon) \\[2mm]
=2\omega \{\bar{v}_2^\epsilon e_1
- \bar{v}_1^\epsilon
e_2+\epsilon[(\bar{v}_2^\epsilon\nabla\bar{\eta}^\epsilon_1-
\bar{v}_1^\epsilon\nabla\bar{\eta}_2^\epsilon)+\bar{B}^\epsilon(\bar{v}_1^\epsilon
 e_2-\bar{v}^\epsilon_2 e_1)]
+\epsilon^2\bar{B}^\epsilon(\bar{v}_1^\epsilon\nabla\bar{\eta}^\epsilon_2-\bar{v}_2^\epsilon\nabla\bar{\eta}^\epsilon_1)\}
  \end{array}
\right.
  \end{eqnarray}
with boundary conditions
\begin{eqnarray*}  %%\label{4.1.7}
\bar{v}^\epsilon_-(t,x',-m)\cdot e_3=\bar{v}^\epsilon_+(t,x',l)\cdot
e_3=0,
\end{eqnarray*}
where
\begin{eqnarray*}  %%\label{4.1.70a}
\begin{array}{ll}
\bar{B}^\epsilon:=(I-(I+\epsilon D\bar{\eta}_\epsilon^T)^{-1})/\epsilon
\end{array}
\end{eqnarray*} and
\begin{eqnarray*}   %%\label{4.1.70b}
\begin{aligned}
\bar{\mathfrak{R}}^\epsilon(t,x)&=\frac1\epsilon
\displaystyle\int_0^{\epsilon\bar{\sigma}^\epsilon(t,x)}
(\epsilon\bar{\sigma}^\epsilon(t,x)-z)h''(\rho_0(x)+z)dz\\[2mm]
&=\epsilon\displaystyle\int_0^{\bar{\sigma}^\epsilon(t,x)}
(\bar{\sigma}^\epsilon(t,x)-z)h''(\rho_0(x)+\epsilon z)dz.
\end{aligned}
\end{eqnarray*}

According to the bound (\ref{4.1.14}) and the sequential weak-$*$ compactness,
we have that up to the extraction of a subsequence (which we still
denote using only $\epsilon$),
\begin{eqnarray*}
 (\bar{\eta}^\epsilon, \bar{v}^\epsilon, \bar{\sigma}^\epsilon)\rightharpoonup
 (\bar{\eta}^0, \bar{v}^0, \bar{\sigma}^0)
 \ {\mbox{weakly-}*\mbox{ in }}  L^\infty((0,t_0),H^3(\Omega)).
 \end{eqnarray*}
 By lower semicontinuity we know that
 \begin{eqnarray}\label{4.17a}
 \sup\limits_{0\leq t<t_0}\|(\bar{\eta}^0,\bar{v}^0, \bar{\sigma}^0)(t)\|_{H^3}\leq 1.
 \end{eqnarray}
Note that by construction $(\bar{\eta}^\epsilon, \bar{v}^\epsilon,
\bar{\sigma}^\epsilon)(0)=(\bar{\eta}, \bar{v}, \bar{\sigma})(0)$.
Our goal is to show that  $(\bar{\eta}^0,\bar{v}^0, \bar{\sigma}^0)$
satisfies the linearized equations (\ref{1.2.14}).

First of all, we give some estimates on $\bar{B}^\epsilon$ and
${\mathfrak{R}}^\epsilon$. We may further assume that $\epsilon$ is
sufficiently small so that
\begin{eqnarray}\label{4.17b}
\sup\limits_{0\leq t<t_0}\|\epsilon
\bar{\sigma}^\epsilon(t)\|_{L^\infty}<\frac12\inf\limits_{-m\leq
x_3\leq l}\rho_0(x_3)
\end{eqnarray}
and $\epsilon<1/(2K_1)$, where $K_1>0$ is the best constant in the inequality
$\|FB\|_{H^4}\leq K_1\|F||_{H^4}\|B\|_{H^4}$
 for $3\times3$ matrix-valued functions $F$, $B$. The former condition implies
 that $\rho_0+\epsilon \bar{\sigma}^\epsilon$ is
 bounded from above and below by positive quantities,
 whereas the latter guarantees that $\bar{B}^\epsilon$
is well-defined and uniformly bounded in $L^\infty(0,t_0;H^2(\Omega))$ since
\begin{equation*}  %%%\label{4.1.15}
\begin{aligned}
\displaystyle\|\bar{B}^\epsilon\|_{H^2}&=\left\|\sum_{n=1}^\infty(-\epsilon)^{n-1}
(D\bar{\eta}^\epsilon)^n\right\|_{H^2}\leq \sum_{n=1}^\infty \epsilon^{n-1}
\|(D\bar{\eta}^\epsilon)^n\|_{H^2}\\[3mm]
&\leq \sum_{n=1}^\infty (\epsilon K_1)^{n-1}\|D\bar{\eta}^\epsilon\|_{H^2}^n\leq\sum_{n=1}^\infty
\frac{1}{2^{n-1}}\|\bar{\eta}^\epsilon\|^n_{H^3}<\sum_{n=1}^\infty\frac{1}{2^{n-1}}=2.
\end{aligned}
\end{equation*}

A straightforward calculation, recalling (\ref{4.17a}) and (\ref{4.17b}), shows that
\begin{eqnarray*}  %%\label{4.1.16}
\sup\limits_{0\leq
t<t_0}\|\bar{\mathfrak{{R}}}^\epsilon(t)\|_{H^3}\leq \epsilon K_3\qquad
\mbox{for some constant }K_3>0.
\end{eqnarray*}

Finally, since ${\rm{Id}}+\epsilon\bar{\eta}^\epsilon$ is invertible, we may define
 $\bar{\zeta}^\epsilon$ via $({\rm{Id}}+\epsilon\bar{\eta}^\epsilon)^{-1}={\rm{Id}}-\epsilon\bar{\zeta}^\epsilon$,
 which implies that $\bar{\zeta}^\epsilon=\bar{\eta}^\epsilon\circ({\rm{Id}}-\epsilon\bar{\zeta}^\epsilon)$.
 The slip map $S^\epsilon_-: {\mathbb{R}}^2\times{\mathbb{R}}^+\rightarrow {\mathbb{R}}^2\times\{0\}$
 is then given by
\begin{eqnarray*} %%\label{4.1.17}
S_-^\epsilon={\rm{Id}}_{{\mathbb{R}}^2}+\epsilon\bar{\eta}^\epsilon_+-\epsilon\bar{\zeta}^\epsilon_+\circ({\rm{Id}}_{{\mathbb{R}}^2}
+\epsilon\bar{\eta}^\epsilon_+).
\end{eqnarray*}
The bounds on $\bar{\eta}^\epsilon$ and the equation satisfied by $\bar{\zeta}^\epsilon$ thus imply that
\begin{eqnarray*} %%\label{4.1.17c}
\sup\limits_{0\leq t<t_0}\|S^\epsilon_-(t)-{\rm{Id}}_{{\mathbb{R}}^2}\|_{L^\infty}\leq
2\epsilon\sup_{0\leq t<t_0}\|\bar{\eta}^\epsilon(t)\|_{L^\infty}
\leq 2\epsilon K_2\sup\limits_{0\leq
t<t_0}\|\bar{\eta}^\epsilon(t)\|_{H^3}<2\epsilon K_2,
\end{eqnarray*}
where $K_2>0$ is the embedding constant for the trace map
$H^3(\Omega)\hookrightarrow L^\infty({\mathbb{R}}^2\times\{0\})$.
This bound allows us to define the normalized slip map
$\bar{S}^\epsilon_-:=(S^\epsilon_--{\rm{Id}}_{{\mathbb{R}}^2})/\epsilon$
as a well-defined and uniformly bounded function in
$L^\infty(0,t_0;L^\infty({\mathbb{R}}^2\times\{0\})).$

 We can now parlay the bounds on $\bar{\eta}^\epsilon$, $\bar{v}^\epsilon$, $\bar{\sigma}^\epsilon$, $\bar{B}^\epsilon$, $\bar{\mathfrak{R}}^\epsilon$
  into corresponding bounds on $\partial_t\bar{\eta}^\epsilon$, $\partial_t\bar{v}^\epsilon$, $\partial_t\bar{\sigma}^\epsilon$,
  and some convergence results. Making use of (\ref{4.141}) and the above bounds of the rescaled functions,
 we infer that
 \begin{eqnarray}
 &&\label{4.1.18}
 \sup\limits_{0\leq t<t_0}\|\partial_t\bar{\eta}^\epsilon(t)\|_{H^3}=\sup\limits_{0\leq t<t_0}\|\bar{v}^\epsilon(t)\|_{H^3}\leq 1,
\\
&& \label{4.1.19}
  \lim\limits_{\epsilon\rightarrow 0}\sup\limits_{0\leq t<t_0}\|\partial_t\bar{\sigma}^\epsilon(t)+\rho_0 {\mathrm{div}}
  \bar{v}^\epsilon(t)\|_{H^2}=0,
\\ && \label{4.1.20} \sup\limits_{0\leq
t<t_0}\|\partial_t\bar{\sigma}^\epsilon(t)\|_{H^2}< K_4,
\\ &&\label{4.1.21}
  \lim\limits_{\epsilon\rightarrow 0}\sup\limits_{0\leq t<t_0}\|\partial_t\bar{v}^\epsilon(t)+\nabla(h'(\rho_0))\bar{\sigma}^\epsilon(t)
  +g e_3\cdot\bar{\eta}^\epsilon(t)-2\omega \bar{v}^\epsilon_2 e_1+2\omega \bar{v}^\epsilon_1
  e_2\|_{H^2}=0,
\\ &&\label{4.1.22} \sup\limits_{0\leq
t<t_0}\|\partial_t\bar{v}^\epsilon(t)\|_{H^2}< K_5.
\end{eqnarray}

We now turn to some convergence results for the jump conditions. We begin with the second equation in (\ref{4.1.5}), which we expand using
the mean-value theorem to get
\begin{eqnarray}\label{4.1.23}
\begin{array}{ll}
&p_+(\rho_0^+)+p'_+(\alpha_+^\epsilon\rho_0^++(1-\alpha^\epsilon_+)\epsilon\bar{\sigma}^\epsilon_+)\epsilon \bar{\sigma}_+^\epsilon\\[2mm]
& =p_-(\rho_0^-)+p'_-(\alpha^\epsilon_-\rho_0^-+(1-\alpha^\epsilon_-)\epsilon\bar{\sigma}^\epsilon_-\circ({\rm{Id}}_{\mathbb{R}^2}+
\epsilon\bar{S}^\epsilon_-))\epsilon\bar{\sigma}_-^\epsilon\circ(\mathrm{Id}_{\mathbb{R}^2}+\epsilon
\bar{S}_-^\epsilon)
\end{array}
\end{eqnarray}
for functions $\alpha^\epsilon_{\pm}: {\mathbb{R}}^+\times
{\mathbb{R}}^2\rightarrow [0,1].$ From the above bounds on
$\bar{\sigma}^\epsilon$ and $\bar{S}^\epsilon$ we get that
 \begin{eqnarray}\label{4.1.24}
 &&\lim\limits_{\epsilon\rightarrow 0}\sup\limits_{0\leq t<t_0}\|
 p'_+(\alpha_+^\epsilon\rho_0^++(1-\alpha^\epsilon_+)\epsilon
 \bar{\sigma}_+^\epsilon)-p'_+(\rho_0^+)\|_{L^\infty}=0,
\\
&& \label{4.1.25}
 \lim\limits_{\epsilon\rightarrow 0}\sup\limits_{0\leq t<t_0}\|p'_-(\alpha_-^\epsilon\rho_0^-+(1-\alpha^\epsilon_-)
 \epsilon\bar{\sigma}_-^\epsilon\circ(\mathrm{Id}_{\mathbb{R}^2}+\epsilon\bar{S}_-^\epsilon))-p'_-(\rho_0^-)\|_{L^\infty}=0
 \end{eqnarray}
and
\begin{eqnarray}\label{4.1.26}
 \lim\limits_{\epsilon\rightarrow 0}\sup\limits_{0\leq t<t_0}\|\bar{\sigma}^\epsilon_-\circ({\rm{Id}}_{R^2}
 +\epsilon\bar{S}^\epsilon_-)-\bar{\sigma}_-^\epsilon\|_{L^\infty}
 \leq\sup\limits_{0\leq t<t_0}\|\nabla\bar{\sigma}^\epsilon(t)\|_{L^\infty}
 \sup\limits_{0\leq t<t_0}\|\epsilon\bar{S}^\epsilon(t)\|_{L^\infty}=0.
 \end{eqnarray}
Since $p_+(\rho_0^+)=p_-(\rho_0^-)$, we may eliminate these terms from the equation (\ref{4.1.23})
and divide both sides by $\epsilon$;
then employing (\ref{4.1.24})--(\ref{4.1.26}), we deduce that
\begin{eqnarray}\label{4.1.27}
\lim\limits_{\epsilon\rightarrow 0}\sup\limits_{0\leq t<t_0}\|p'_+(\rho_0^+)
\bar{\sigma}^\epsilon_+(t)-p'_-(\rho_0^-)\bar{\sigma}^\epsilon_-(t)\|_{L^\infty}=0.
\end{eqnarray}

For the second equation in (\ref{4.1.5}) we first write the normal at the interface as $n^\epsilon=N^\epsilon/|N^\epsilon|$ with
\begin{eqnarray*}  %%\label{4.1.28}
\begin{array}{ll}
N^\epsilon&=(e_1+\epsilon\partial_{x_1}\bar{\eta}^\epsilon_+)\times(e_2+\epsilon\partial_{x_2}\bar{\eta}^\epsilon_+)\\[2mm]
&=e_3+\epsilon(e_1\times\partial_{x_2}\bar{\eta}^\epsilon_++\partial_{x_1}\bar{\eta}^\epsilon_+\times
e_2)+\epsilon^2(\partial_{x_1}\bar{\eta}^\epsilon_+
\times\partial_{x_2}\bar{\eta}^\epsilon_+)=:e_3+\epsilon\bar{N}^\epsilon.
\end{array}
\end{eqnarray*}
As $\epsilon\rightarrow 0$ we have $|N^\epsilon|>0$, so we may rewrite the first equation in (\ref{4.1.5}) as
\begin{eqnarray*}  %%\label{4.1.29}
(\bar{v}^\epsilon_+-\bar{v}^\epsilon\circ({\rm{Id}}_{\mathbb{R}^2}+\epsilon\bar{S}^\epsilon))\cdot
(e_3+\epsilon\bar{N}^\epsilon)=0.
\end{eqnarray*}
Clearly $\sup\limits_{0\leq t\leq t_0}\|\bar{N}^\epsilon(t)\|_{L^\infty}$ is bounded uniformly and
\begin{eqnarray*}  %%\label{4.1.30}
 \lim\limits_{\epsilon\rightarrow 0}\sup\limits_{0\leq t<t_0}\|\bar{v}^\epsilon
 \circ({\rm{Id}}_{\mathbb{R}^2}+\epsilon\bar{S}^\epsilon)-\bar{v}_-^\epsilon\|_{L^\infty}
 \leq\sup\limits_{0\leq t<t_0}\|D\bar{v}^\epsilon(t)\|_{L^\infty}
 \sup\limits_{0\leq t<t_0}\|\epsilon\bar{S}^\epsilon(t)\|_{L^\infty}=0,
\end{eqnarray*}
from which we find that
\begin{eqnarray}\label{4.1.31}
\lim\limits_{\epsilon\rightarrow 0}\sup\limits_{0\leq t<t_0}\|e_3
\cdot(\bar{v}^\epsilon_+(t)-\bar{v}^\epsilon_-(t))\|_{L^\infty}=0.
\end{eqnarray}

 According to (\ref{4.1.18}), (\ref{4.1.20}) and (\ref{4.1.22}), by Lions-Aubin Lemma,
 we see that the set $\{(\bar{\eta}^\epsilon,\bar{\sigma}^\epsilon,\bar{v}^\epsilon)\}$
 is strongly pre-compact in the spaces $L^\infty(0,t_0;H^{\frac{11}{4}})$. Therefore,
  \begin{eqnarray*}  %%\label{4.1.39}
  (\bar{\eta}^\epsilon, \bar{\sigma}^\epsilon, \bar{v}^\epsilon)\rightarrow
  (\bar{\eta}^0,\bar{\sigma}^0,\bar{v}^0)\  {\rm{strongly\  in}}\
   L^\infty(0,t_0;H^{\frac{11}{4}}).
  \end{eqnarray*}

  This strong convergence, together with the convergence results (\ref{4.1.19}) and
  (\ref{4.1.21}) and the equation
  $\partial_t\bar{\eta}^\epsilon=\bar{v}^\epsilon$, implies that
  \begin{eqnarray*}   %%% \label{4.1.40}
  \left(\partial_t\bar{\eta}^\epsilon,\partial_t\bar{v}^\epsilon,\partial_t\bar{\sigma}^\epsilon\right)
  \rightarrow \left(\partial_t\bar{\eta}^0,\partial_t\bar{v}^0,\partial_t\bar{\sigma}^0\right) {\rm{strongly\ in}}
  \ L^\infty(0,t_0;H^\frac74(\Omega)),
  \end{eqnarray*}
  and that
  \begin{eqnarray*}   %% \label{4.1.41}
  \left\{
  \begin{array}{ll}
  \partial_t\bar{\eta}^0=\bar{v}^0,\\[2mm]
  \partial_t\bar{\sigma}^0+\rho_0{\mathrm{div}} \bar{v}^0=0\\[2mm]
  \partial_t\bar{v}^0+\nabla(h'(\rho_0)\bar{\sigma}^0+g e_3\cdot\bar{\eta}^0)=2\rho_0\omega\bar{v}^0_2 e_1-2\rho_0\omega \bar{v}^0_1 e_2.
  \end{array}
  \right.
  \end{eqnarray*}
  We may pass to the limit in the initial conditions
  $(\bar{\eta}^\epsilon,\bar{v}^\epsilon,\bar{\sigma}^\epsilon)(0)=(\bar{\eta},\bar{v},\bar{\sigma})(0)$
  to find that
  \begin{eqnarray}\label{4.1.42}
  (\bar{\eta}^0,\bar{v}^0,\bar{\sigma}^0)(0)=(\bar{\eta},\bar{v},  \bar{\sigma})(0)
  \end{eqnarray}
  as well.

  We now derive the jump and boundary conditions for the limiting functions.
  The index $11/4$ is sufficiently large to give the
  $L^\infty(0,t_0;L^\infty)$-convergence of $(\bar{\eta}^\epsilon$, $\bar{v}^\epsilon$,
  $\bar{\sigma}^\epsilon)$ when restricted to $\{x_3=0\}$, $\{x_3=-m\}$, and $\{x_3=l\}$,
  i.e., the interface and the lower and upper boundaries.
  Combining this with (\ref{4.1.27}) and (\ref{4.1.31}), we deduce that
 %% \begin{eqnarray*}\label{4.1.43}
 $$ p'_+(\rho_0^+)\bar{\sigma}^0_+=p'_-(\rho_0^-)\bar{\sigma}^0_- \ {\rm{on}}\ \{x_3=0\},
 %% \end{eqnarray*}   \begin{eqnarray*}\label{4.1.44}
 \qquad (\bar{v}_+^0-\bar{v}_-^0)\cdot e_3=0 \ {\rm{on}}\ \{x_3=0\}, $$
%%  \end{eqnarray*}
  and
  \begin{eqnarray*}   %% \label{4.1.45}
  \bar{v}^0_+\cdot e_3=0 \ {\rm{on}}\ \{x_3=l\}, \qquad \bar{v}^0_-\cdot e_3=0 \ {\rm{on}}\ \{x_3=-m\}.
  \end{eqnarray*}

Henceforth, $(\bar{\eta}^0,\bar{v}^0,\bar{\sigma}^0)$ is a solution
to (\ref{1.2.14}) along with the corresponding jump and boundary
conditions on $(0,t_0)\times\Omega$ which satisfies the initial
condition (\ref{4.1.42}). Thus, in view of the uniqueness of the
linearized equations Theorem 3.1, we have
\begin{eqnarray*}  %%%\label{4.1.40}
(\bar{\eta}^0,\bar{v}^0,\bar{\sigma}^0)=(\bar{\eta},\bar{v},\bar{\sigma})
\mbox{ on } [0,t_0)\times \Omega.
\end{eqnarray*}
Hence we may chain together inequalities (\ref{4.1.12a}) and (\ref{4.17a}) to get
\begin{eqnarray*}   %%\label{4.1.41}
2<\sup\limits_{t_0/2\leq t<t_0}\|(\bar{\eta}^0,\bar{v}^0,\bar{\sigma}^0)(t)\|_{H^3}
\leq\sup\limits_{0 \leq t<t_0}\|(\bar{\eta}^0,\bar{v}^0,\bar{\sigma}^0)(t)\|_{H^3}\leq 1,
\end{eqnarray*}
which is a contradiction. Therefore, the perturbed problem does not have property
$EE(k)$ for any $k\geq 3. $\hfill $\Box$
\par \quad \\
 {\bf{Acknowledgement}}\quad The authors are grateful to the
referees for their helpful suggestions which improved the
presentation of this paper.

\end{document}